\documentclass[11pt]{article}
\usepackage{amsmath}
\usepackage{amssymb}
\usepackage{amsfonts} 
\usepackage{amstext} 
\usepackage{amsthm} 
\usepackage{euscript} 
\usepackage{eufrak} 

\newcommand{\Z}{\ensuremath{\mathbb{Z}}}
\newcommand{\N}{\ensuremath{\mathbb{N}}}
\newcommand{\Q}{\ensuremath{\mathbb{Q}}}
\newcommand{\R}{\ensuremath{\mathbb{R}}} 
\newcommand{\C}{\ensuremath{\mathbb{C}}}

\newcommand{\cC}{\ensuremath{\mathcal{C}}}
 
\newcommand{\cE}{\ensuremath{\mathcal{E}}}

\newcommand{\cH}{\ensuremath{\mathcal{H}}}

\newcommand{\cL}{\ensuremath{\mathcal{L}}}
\newcommand{\cM}{\ensuremath{\mathcal{M}}} 
\newcommand{\cN}{\ensuremath{\mathcal{N}}} 
\newcommand{\cP}{\ensuremath{\mathcal{P}}} 
 
\newcommand{\cS}{\ensuremath{\mathcal{S}}}
\newcommand{\cA}{\ensuremath{\mathcal{A}}}

\newcommand{\hh}{\ensuremath{\hat{h}}} 
\newcommand{\he}{\ensuremath{\hat{e}}} 
\newcommand{\hU}{\ensuremath{\hat{\Upsilon}}}

\newcommand{\fC}{\ensuremath{\EuFrak{C}}}
\newcommand{\fT}{\ensuremath{\EuFrak{T}}}
\newcommand{\fB}{\ensuremath{\EuFrak{B}}}

\newcommand{\sS}{\ensuremath{\EuScript{S}}}

\begin{document} 

\begin{center} 
{\Large \bf A Holomorphic 
0-Surgery Model for Open Books with Application 
to Cylindrical Contact Homology} \\ 

\vspace{.3in} 
Mei-Lin Yau  \\  
   \vspace{.2in} 
Department of Mathematics \\  Michigan State University \\ 
East Lansing, MI 48824 \\ 
Email: yau@math.msu.edu
\end{center} 
\vspace{.2in} 
\begin{abstract} 
We give a simple model in $\C^2$ of the 
0-surgery along a fibered knot of a closed 3-manifold 
$M$ to yield a mapping torus $\hat{M}$. This model allows 
explicit relations between pseudoholomorphic curves in 
$\R\times M$ and in $\R\times \hat{M}$. 
We then use it to compute the cylindrical contact homology 
of open books resulting from a positive Dehn twist on 
a torus with boundary. 
\end{abstract}


\newtheorem{cond}{Condition}[section]
\newtheorem{defn}{Definition}[section]
\newtheorem{exam}{Example}[section]
\newtheorem{lem}{Lemma}[section]
\newtheorem{cor}{Corollary}[section]
\newtheorem{main}{Main Theorem}
\newtheorem{theo}{Theorem}[section]
\newtheorem{prop}{Proposition}[section]
\newtheorem{rem}{Remark}[section]
\newtheorem{notn}{Notation}[section] 
\newtheorem{fact}{Fact}[section] 





\section{Introduction} \label{intro} 

In this note all 3-manifolds are closed and orientable,  
all contact structures are coorientable and all surfaces 
are oriented. Given a surface $\Sigma$ a diffeomorphism 
$\phi\in \text{Diff}^+(\Sigma)$ we denote by 
\[ 
\Sigma_\phi:=\frac{\Sigma\times [0,1]}{(\phi(x),0)\sim (x,1)} 
\] 
the mapping torus associated to $\phi$. Throughout this paper we 
assume that a symplectic structure $\omega$ on $\Sigma$ is given, 
$\phi\in\text{Symp}(\Sigma,\omega)$ and 
$\phi=id$ near $\partial\Sigma$.

It is a basic fact in topology that a canonical 0-surgery along 
(every connected component of) the binding of an open book $M$ 
yields a mapping torus $\hat{M}$. Now with the correspondence 
between contact structures and open books established by 
Thurston and Winkelnkemper \cite{TW} and Giroux \cite{G1,GM}, 
it is expected that a "nice" description of the said 0-surgery 
will benefit the study of contact manifolds and symplectic 
manifolds.

For example, Eliashberg \cite{E3} showed that 
this surgery can be done {\em symplectically}, namely the two 
manifolds can be included as the boundary of a symplectic 
cobordism of which the symplectic structure satisfies some 
boundary conditions pertaining to the given open book and mapping 
torus. His result leads to the equivalence between the 
(weakly) semi-symplectic
fillability and the (weakly) symplectic fillability of  
contact manifolds, and provides applications to 
Kronheimer and Mrowka's Property P as well as  
Ozsv\'{a}th and Szab\'{o}'s Heegaard Floer homology theory. 
See also \cite{Et,O}.

Now, both $M$ and $\hat{M}$ have a natural 
{\em symplectization} $\R\times M$ and $\R\times\hat{M}$ on 
which one can define holomorphic curve invariants with similar 
setups. On $\R\times M$ we have {\em contact homology} first 
constructed by Eliashberg and Hofer 
\cite{E2,EGH} to provide Gromov-Floer 
type invariants for contact manifolds. On $\R\times \hat{M}$ 
there is Hutchings and Sullivan's {\em periodic Floer homology} for
symplectic maps \cite{HS}, 
which is a generalization of Seidel's {\em symplectic 
Floer homology} \cite{Se}. 

Here we are interested in the contact homology of 
contact 3-manifolds, 
which is  in general very difficult to compute. 
To provide an access for computing and 
studying contact homology, 
it is then 
desirable to find a model demonstrating that the
0-surgery  can be done {\em holomorphically} -- at least in some
reasonable sense, allowing an explicit correspondence between moduli
of pseudoholomorphic curves and hence, a comparison of 
pseudoholomorphic curve theories 
on the two symplectic manifolds.

Indeed, such a model does exist and is very simple. Assume the 
binding $B$ of $M$ is connected for simplicity. Let $\cN_B$ denote 
a small tubular neighborhood of $B$. Let $\hat{\cN}_B:=\hat{M}
\setminus (M\setminus \cN_B)$. $\hat{\cN}_B$ is a tubular 
neighborhood of an orbit $\hat{e}$ corresponding to an 
elliptic fixed point of the monodromy of the mapping torus $\hat{M}$. 

Let $\C^2$ be 
the standard complex plane. We will prove the following 

\begin{theo} \label{main1} 
There are simultaneous 
holomorphic embeddings of symplectizations $\R\times \cN_B$ and 
$\R\times \hat{\cN}_B$ into $\C^2$ with the intersection 
of the images an open domain in $\{ z_1z_2\neq 0\}$
\end{theo} 

In particular, with appropriate almost complex structures on 
$\R\times M$ and 
$\R\times \hat{M}$ given, Theorem  \ref{main1} implies  

\begin{lem} \label{mainlem} 
Let $C\subset \R\times M$ 
be a pseudoholomorphic cylinder bounding Reeb orbits 
$\gamma_{\pm}$ at $\pm\infty$, with $\gamma_{\pm}$ not equal 
to any multiple of $B$. Suppose that $C$ intersects with $\R\times B$ 
at $s$ points with intersection multiplicities $m_1,...,m_s\in \N$, 
then $C$ lifts, via the canonical 0-surgery along $B$, 
to a $(2+s)$-punctured pseudoholomorphic sphere in $\R\times \hat{M}$. 
Moreover, the extra $i^{th}$ puncture 
converge to the $m_i^{th}$ iterate 
of $\hat{e}$ at $-\infty$. 
\end{lem}

Let $(M,\xi)$ be the contact 3-manifold associated to the 
open book $(\Sigma,\phi)$, where $\Sigma$ is a punctured 
torus and the monodromy $\phi$ is a positive $\sigma$-Dehn twist 
along an embedded  nonseparating circle of $\Sigma$. 
We have $H_1(M,\Z)=\Z_\sigma\oplus \Z$. 
The holomorphic 0-surgery model enables us to relate certain 
holomorphic cylinders in $\R\times M$ to Taubes's 
trice-punctured spheres \cite{Ta} (see also \cite{HS}), and 
by using Bourgeois's Morse-Bott version of contact homology  
We obtain 
the following 

\begin{theo} \label{example} 
The cylindrical contact homology $HC(M,\xi,\Q)$ is 
freely generated by 
\begin{enumerate} 
\item $h^m$ (hyperbolic) with $m\in\N$, $[h^m]=0\in H_1(M,\Z)$, and by 
\item $E_{i,m}$ (elliptic) with $(i,m)\in\Z_\sigma\times\N\setminus \{
(0,1)\}$, 
$[E_{i,m}]=(i,0)\in H_1(M,\Z)$. 
\end{enumerate} 
Moreover, $[h^m]=[E_{0,m}]=0\in H_1(M,\Z)$, and 
their reduced 
Conley-Zehnder indexes are 
\[ 
\bar{\mu}(h^m)=2m-1 \ \ m\in\N; \quad \bar{\mu}(E_{0,m})=2m-2, 
\ \ m\in\N_{\geq 2}. 
\] 
\end{theo} 

Readers are referred to Proposition \ref{bct} and 
Definition \ref{E_i} for the definitions of $
h^m$ and $E_{i,m}$. 

Note that when the positive Dehn twist is simple 
($\sigma=1$), the open 
book is $S^1\times S^2$ with the unique (up to isotopy) 
Stein-fillable contact structure, of which   
the cylindrical contact has been computed in 
\cite{Y2} via {\em subcritical}  
contact handle attaching. However, for $\sigma\geq 2$ 
the contact manifold is Stein-fillable but {\em not} 
subcritical Stein-fillable, nor an $S^1$-bundle over a 
closed surface (the contact homology of $S^1$-bundles 
have been computed, see \cite{EGH}\cite{B}). Our 
results here provide first examples of cylindrical 
contact homology via nontrivial (i.e. monodromy $\neq id$) 
open books. We hope that the holomorphic 0-surgery model 
will lead to 
more examples of (cylindrical) contact homology of open books 
as well as new results in contact topology.

This paper is organized as follows. Section \ref{back} 
consists of some background on cylindrical contact homology, 
contact structures associated 
to open books, 0-surgery and mapping tori. In Section 
\ref{model} we construct in $\C^2$ a holomorphic model of a 
0-surgery and verify Theorem \ref{main1} and Lemma \ref{mainlem}. 
The computation of the cylindrical contact homology of a 
positive Dehn twist is done in 
Section \ref{appl}. 

\vspace{.3in} 
\noindent 
{\bf \Large Acknowledgement:} 

The author is grateful to the referee of an earlier version 
of this paper for critical comments and many valuable 
suggestions.

\section{Background} \label{back}

\subsection{Cylindrical contact homology} 

\noindent 
{\bf Contact forms.} \ 
A 1-form $\alpha\in\Omega^1(M)$ ($\dim M=3$) is said to be {\em contact}
if 
$\alpha\wedge d\alpha$ is nowhere vanishing. The kernel 
$\xi:=\ker \alpha$ is called a {\em contact structure}. We say 
$\alpha$ and hence $\xi$ are {\em positive} if 
$\alpha\wedge d\alpha$ is a volume form of the oriented manifold 
$M$. In this 
paper all contact 1-forms considered are {\em positive}. 

\vspace{.2in} 
\noindent 
{\bf Reeb orbits.} \ 
There associates to $\alpha$ a unique vector 
field $R=R_\alpha$  called {\em Reeb} vector field, which is 
defined by 
\[ 
d\alpha(R,\cdot )=0, \quad \alpha(R)=1. 
\] 

A periodic integral trajectory of $R$ is called a {\em Reeb orbit} 
(of $\alpha$). We call $\gamma$ {\em simple} if $\gamma$
is not a nontrivial multiple cover of another Reeb orbit. 

\begin{notn} 
{\rm 
We denote by $\gamma^m$ the $m^{th}$-iterate of 
a Reeb orbit $\gamma$. 
} 
\end{notn} 

\begin{defn}[action] 
{\rm 
Let $\gamma :[o,\tau ]\to M$ be a Reeb trajectory with
$\dot{\gamma}(t)=R_{\alpha}(\gamma (t))$. Define the {\em action}
${\cal A}(\gamma )$ of $\gamma$ to be the number
\begin{equation} \label{act} 
T ={\cal A}(\gamma ):=\int _{\gamma}\alpha
\end{equation} 
} 
\end{defn} 

The flow $R^t$ of $R$ preserves $\xi$.
Thus the linearized Reeb flow
$R^t_*$ , when restricted on $\gamma$,
defines a path of symplectic maps
\[
\Lambda_\gamma (t)=R^t_*
(\gamma (0)):\xi |_{\gamma (0)}\to \xi |_{\gamma (t)} \ .
\]

\vspace{.2in} 
\noindent 
{\bf Linearized Poincar\'{e} return map.} \ 
When $\gamma$ is a Reeb orbit with action $T$,
$\Lambda_\gamma:=\Lambda_\gamma (T)$ is called the 
{\em linearized Poincar\'{e}
return  map} along $\gamma$. 
\begin{defn} 
{\rm 
A Reeb orbit $\gamma$ is {\em non-degenerate}
if 1 is not an eigenvalue of its linearized 
Poincar\'{e} return map $\Lambda_\gamma$. 
A contact 1-form $\alpha$ is called {\em regular}
if every Reeb orbit of $\alpha$ is non-degenerate. 
} 
\end{defn}   
It is well-known that generic contact 1-forms are regular 
(see \cite{B2}).

\begin{defn}[good orbit]  \label{good} 
{\rm 
A Reeb orbit is said to be {\em bad} 
(see Section 1.2 of \cite{EGH}) if it is an even multiple 
of another Reeb orbit whose 
linearized Poincar\'{e} return map has the 
property that the total multiplicity of its eigenvalues 
from the interval $(-1,0)$ is odd. A Reeb orbit is {\em good} 
if it is not bad. 
} 
\end{defn} 

\begin{notn} 
{\rm 
We denote by $\cP_{\alpha}$ the set of all {\em good} Reeb orbits 
of $\alpha$. } 
\end{notn} 

\vspace{.2in} 
\noindent 
{\bf A mod 2 index.} \  
Assume that $\alpha$ is regular. Then 
for {\em any} Reeb orbit $\gamma$, a $\Z_2$-index 
$\bar{\mu}(\gamma, \Z_2)$ is defined: 
\begin{equation}  \label{z2ind}
\bar{\mu}(\gamma;\Z_2)=\begin{cases} 
0 \quad \text{if $\gamma$ is {\em even}, i.e., if $\det (\Lambda_\gamma
-Id)>0$};  \\ 
1 \quad \text{if $\gamma$ is {\em odd}, i.e., if $\det (\Lambda_\gamma
-Id)<0$}. 
\end{cases} 
\end{equation} 

\vspace{.2in} 
\noindent 
{\bf $\bar{\mu}$-index.} \ 
Since $d\alpha |_\xi$ is a symplectic 
2-form, the first Chern class $c_1(\xi)\in H^2(M,\Z)$ 
is defined. 
For the sake of simplicity we assume in this paper that 
$c_1(\xi)=0$ on $H_2(M,\Z)$. Then the Conley-Zehnder index 
$\mu({\gamma})$ of a {\em homologously trivial} Reeb orbit $\gamma$ 
is well-defined (\cite{SZ}\cite{RS}). 

\begin{defn}
{\rm 
The 
{\em reduced} Conley-Zehnder index is defined to be 
\begin{equation}  \label{reduced-mu} 
\bar{\mu}(\gamma):=\mu(\gamma)-1 \quad \text{ if }\dim M=3.  
\end{equation} 

When $\gamma$ is homologously trivial, 
\begin{equation} \label{mumu} 
\bar{\mu}(\gamma)\equiv \bar{\mu}(\gamma, \Z_2) \quad 
\mod 2. 
\end{equation} 
}
\end{defn} 

\vspace{.2in} 
\noindent 
{\bf Almost complex structures.} \ 
Following Eliashberg, Givental and Hofer \cite{E2}\cite{EGH} 
one can define Gormov-Floer type invariant called {\em contact
homology}  for $(M,\xi)$ by counting in the symplectic manifold 
(called the {\em symplectization} of $(M,\alpha)$) 
$(\R\times M,d(e^t\alpha))$ pseudoholomorphic curves bounding 
Reeb orbits at $\pm\infty$. The almost complex structure $J$ 
involved is $\alpha$-{\em admissible}, i.e., 
\begin{enumerate} 
\item $J$ preserves $\xi$;  
\item $J|_\xi$ is $d\alpha$-compatible, i.e., $d\alpha (v,Jv)>0$ 
for all $0\neq v\in \xi$ and $d\alpha (v,w)=d\alpha(Jv,Jw)$ for 
all $v,w\in \xi$; 
\item $J(\partial_t)=R_\alpha$, $J(R_\alpha)=-\partial_t$. 
\end{enumerate} 

\vspace{.2in} 
\noindent
{\bf Pseudoholomorphic cylinders and planes.} \ 
Fix a pair $(\alpha ,J)$ with $\alpha$ regular and $J$ an 
$\alpha$-admissible almost complex structure on $\R\times M$. 
Given two good Reeb orbits $\gamma_-$ and $\gamma_+$  
we denote by $\cM(\gamma _-,\gamma _+)$ 
the moduli space of maps $(\tilde{u},j)$ where 
\begin{enumerate} 
\item $j$ is an almost complex structure on 
$S^2$ (here we identify $S^2$; 
\item let $\dot{S}^2:=S^2\setminus \{ 0,\infty\}$, then 
$\tilde{u}=(a,u):(\dot{S}^2,j)\to 
(\R\times M,J)$ is a proper map and is $(j,J)$-holomorphic, 
i.e., $\tilde{u}$ satisfies
$d\tilde{u}\circ j= J\circ d\tilde{u}$; 
\item $\tilde{u}$ is asymptotically cylindrical over $\gamma _-$ 
at the negative end of $\R\times M$ at the puncture 
$0\in S^2$; and $\tilde{u}$ is asymptotically cylindrical 
over $\gamma _+$ at the positive end of $\mathbb{R}\times M$ 
at the puncture $\infty\in S^2$;  
\item $(\tilde{u},j)\thicksim (\tilde{v},j')$ if there is a 
diffeomorphism $f:\dot{S}^2\to \dot{S}^2$ such that 
$\tilde{v}\circ f=\tilde{u}$, $f_*j=j'$, 
and $f$ fixes all punctures. 
\end{enumerate} 

For a {\em contractible} Reeb orbit $\gamma$ 
The moduli space $\cM(\gamma)$ of pseudoholomorphic planes 
bounding $\gamma$ at $\infty$ 
is defined in a similar fashion.

For generic $\alpha$-admissible $J$, $\cM(\gamma_-,\gamma_+)$ 
and $\cM(\gamma)$, 
if not empty, are smooth manifolds on which $\R$ 
acts freely by translation. If both 
$\gamma_{\pm}$ are homologously trivial then (see \cite{EGH}
Proposition 1.7.1) 
\begin{equation} \label{dim+-}  
\dim\cM (\gamma_-,\gamma_+)=\bar{\mu}(\gamma_+)-\bar{\mu}
(\gamma_-).  
\end{equation} 
In particular, if $\dim\cM (\gamma_-,\gamma_+)=1$ then 
$\cM (\gamma_-,\gamma_+)/\R$ is 
compact $0$-dimensional, 
hence a finite number of points. 

For contractible $\gamma$ we have 
\begin{equation} \label{dimdisc} 
\dim \cM(\gamma)=\bar{\mu}(\gamma). 
\end{equation}

\vspace{.2in} 
\noindent 
{\bf Energy.} \ 
If $\tilde{u}=(a,u)\in\cM(\gamma _-,\gamma _+)$ (or 
$\cM(\gamma)$) then 
$u^*d\alpha \geq 0$ pointwise, and vanishes at most at finitely 
many points. We define the {\em contact energy} $E(\tilde{u})$ 
of $\tilde{u}$ to be 
\[ 
E(\tilde{u}):=\int _{u(\dot{S}^2)}d\alpha =\int _{\gamma_+}
\alpha -\int_{\gamma_-}\alpha=\cA(\gamma_+)-\cA(\gamma_-)\geq 0. 
\] 
Note that 
$E(\tilde{u})=0$ iff $\gamma _-=\gamma _+$, and in this case 
the moduli space consists of a single element 
$\R\times\gamma _+$. 

For $\tilde{u}=(a,u)\in\cM(\gamma)$ the contact energy is 
defined similarly: 
\[ 
E(\tilde{u}):=\int _{u(\C)}d\alpha =\int _{\gamma}
\alpha =\cA(\gamma) >0. 
\] 

\vspace{.2in} 
\noindent
{\bf Contact complex.} \ 
The contact complex 
$\cC(\alpha)$ is the free module over $\Q$ generated by 
all elements of $\cP_\alpha$ the set of all {\em good} Reeb orbits.

\vspace{.2in} 
\noindent
{\bf Boundary operator $\partial$.} \ 
For a Reeb orbit $\gamma$ we 
denote by $\kappa_\gamma$ its multiplicity. Similarly we 
denote by $\kappa_C$ the multiplicity of a pseudoholomorphic 
curve $C$ in $\R\times M$. 

The {\em boundary operator} $\partial$ 
of the contact complex $\cC(\alpha)$ is defined by 
(see \cite{E3}\cite{B2} but for a different coefficient ring) 
\begin{align}  
\partial \gamma &:=
{\sum_{\gamma'\in\cP_\alpha}}\langle \partial
\gamma,\gamma'\rangle \gamma' \\ 
 \langle \partial
\gamma,\gamma'\rangle &:= \kappa_\gamma 
\underset{\dim
\cM(\gamma',\gamma)=1}{\sum_{C\in\cM(\gamma',\gamma)/\R}} 
\frac{\pm 1}{\kappa_C}   \label{coeff} 
\end{align}   
The $\pm$ sign 
in (\ref{coeff}) depends on the orientation of 
$C\in\cM(\gamma',\gamma)/\R$ (see Section \ref{orien}).

\begin{defn} 
{\rm 
Suppose that $\partial^2=0$. Then the {\em cylindrical contact
homology}  of $(M,\xi,\alpha, J)$ is defined to be
$HC(M,\xi,\alpha,J):=\ker{\partial}/{\rm im}\partial$. 
} 
\end{defn}

\begin{rem} \label{d} 
{\rm  In contact homology \cite{EGH} one defines the boundary 
operator $d=\sum_{i=0}^{\infty}d_i$ by counting 1-dimensional 
moduli of holomorphic spheres with one positive puncture and 
arbitrary number of negative punctures. The summand $d_i$ counts 
the number of holomorphic spheres with $i$ negative punctures. 
In particular $d_1=\partial$. 
Since $d^2=0$ (see \cite{EGH}) we have $d_1^2+d_0d_2=0$. 
Thus $\partial^2=0$ if $d_0=0$, which is the case for the 
examples that we will compute in Section \ref{appl}. 
} 
\end{rem} 

\vspace{.2in} 
\noindent 
{\bf Contractible subcomplex.} \ 
Let $\cC^o(\alpha)$ denote the subcomplex generated by all good 
{\em contractible} Reeb orbits. Recall that we assume 
$c_1(\xi)=0$ on $H_2(M,\Z)$, hence $\bar{\mu}$-index is well-defined 
for all contractible Reeb orbits. 
Thus $\cC^o(\alpha)$ is graded by $\bar{\mu}$. 
We denote by $\cC_k^o(\alpha)\subset \cC^o(\alpha)$ the subcomplex 
generated by all elements of $\cC^o(\alpha)$ with $\bar{\mu}=k$.

\begin{theo}[see \cite{U2}\cite{EGH}]  \label{inv} 
Assume that $\partial^2=0$. Suppose that $\cC^o_k(\alpha)=0$ for 
$k=0,-1$, then the cylindrical contact homology $HC(M,\xi)
:=HC(M,\xi,\alpha,J)$ is 
independent of the contact form $\alpha$, the almost complex 
structure $J$; it depends only on the isotopy class of the 
contact structure $\xi$.  
\end{theo}

\subsection{Open book, contact structure and mapping torus}

\noindent 
{\bf Open book.} \  
The pair $(\Sigma, \phi)$ is said to be an 
{\em open book} representation of a 3-manifold $M$ if 
$M$ can be expressed as 
\[ 
M=\Sigma_\phi\cup_{id} (B\times D^2) 
\] 
where 
\begin{itemize} 
\item $\phi\in\text{Diff}^+(\Sigma,\partial\Sigma)$, $\phi=id$ near 
$\partial\Sigma$ is the {\em monodromy}, and 
\item $B\cong \partial\Sigma$ is called the {\em binding} of the 
open book. 
\end{itemize} 
The complement $M\setminus B$ fibers over $S^1$, the fibers are 
called {\em pages} and are diffeomorphic to $\Sigma$. 

\vspace{.2in} 
\noindent 
{\bf Positive stabilization.} \ 
An open book $(\Sigma',\phi')$ is called a {\em positive
stabilization}   of $(\Sigma, \phi)$ if 
\begin{itemize} 
\item 
$\Sigma'$ is obtained by gluing 
to $\partial \Sigma$ the 
boundaries $\{ \pm 1\}\times [-\epsilon,\epsilon]$ of 
a strip (a 2-dimensional handle of index 1) $[-1,1]\times 
[-\epsilon,\epsilon]$  (so the Euler number of 
the page is decreased by 1, i.e., $\chi(\Sigma')=\chi(\Sigma)-1$), 
\item 
$\phi':=\phi\circ\tau_\Gamma$; where $\Gamma\subset \Sigma'$ is an 
embedded circle intersecting with the cocore $\{ 0\}\times [-\epsilon, 
\epsilon]$ at a single point, and $\tau_\Gamma$ is the 
{\em positive} Dehn twist along $\Gamma$.  
\end{itemize} 

Two open books are said to be {\em equivalent up to positive 
stabilizations} if they become equal after applied with 
finitely many positive stabilizations. 

With positive stabilizations we may assume that {\em $B$ is 
 connected} \cite{G1}.

\vspace{.2in} 
\noindent
{\bf Associated contact structure.} \ 
The work of Thurston and Winkelnkemper \cite{TW} and 
Giroux \cite{G1} established the following important bijection: 
\[ 
\frac{\text{contact structures on $M$}}{\text{contact isotopies}}
\leftrightarrow \frac{\text{open books of $M$}}{\text{positive 
stabilizations}}. 
\] 

Here we sketch the construction of a contact 1-form 
associated to an open book (see \cite{G1}). 
For an open book $(\Sigma,\phi)$ with a connected 
binding $B$, we fix an area form 
$\omega=d\beta$ on $\Sigma$. Isotope $\phi$ 
if necessary we may assume that 
$\phi^*\omega =\omega$. 
Then associate to $(\Sigma,\phi)$ a contact 
1-form $\alpha$ such that 
\begin{itemize} 
\item $d\alpha$ 
restricts to an area form on every page,  
$d\alpha=\omega$ on each fiber of $\Sigma_\phi$, and 
\item $\alpha=dp+r^2dt$ near the binding $B$, 
where $p$ parametrizes $B\cong S^1$, $(r,t)$ are the polar coordinates 
of the $D^2$ factor of $B\times D^2$. 
\end{itemize}   
In particular, on $\Sigma_\phi$ we can define $\alpha$ to be 
\[ 
\alpha := (1-t)\beta +t\phi^*\beta +Kdt. 
\] 
Then $\alpha$ 
is contact on $\Sigma_\phi$ provided that $K$ is a large 
enough constant. 

\begin{rem} \label{Rorbit} 
{\rm 
If $\alpha$ is associated to an open book $(\Sigma,\phi)$ then 
the Reeb orbits of $\alpha$ 
correspond to the the periodic points of 
some $\phi'\in\text{Symp}(\Sigma,\omega)$  isotopic to 
$\phi$, 
and positive multiples of $B$.  
} 
\end{rem}

\vspace{.2in} 
\noindent
{\bf Mapping torus from 0-surgery.} \ 
A canonical 
0-surgery along the binding $B$ of an open book $(\Sigma,\phi)$ of 
a 3-manifold $M$ 
yields the mapping torus $\hat{M}=\hat{\Sigma}_{\hat{\phi}}$ 
together with a special section $\hat{e}=\{ z\}\times S^1$ 
of the fibration 
\[ 
\hat{\Sigma}\to \hat{M}\overset{\hat{\pi}}{\to}S^1, 
\] 
where 
\begin{itemize} 
\item $\hat{\Sigma}=\Sigma\cup D^2$ is a closed Riemann surface, 
\item $z\in D^2\subset \hat{\Sigma}$, $\hat{M}\setminus \hat{e} 
\overset{\text{diffeo}}{\cong}\Sigma_\phi$, 
\item 
$\hat{\phi}\in\text{Symp}(\hat{\Sigma},\hat{\omega})$ satisfies  
$\hat{\omega}|_\Sigma=\omega$, 
$\hat{\phi}|_{\Sigma}=\phi$, $\hat{\phi}|_{D^2}=id$ 
(so $\hat{\phi}(z)=z$).  
\end{itemize} 

\begin{rem} 
{\rm 
The point $z\in\text{Fix}(\hat{\phi})$ and hence the 
section $\hat{e}=\{ z\}\times S^1$ remember the 0-surgery, so that a 
canonical 0-surgery along $\hat{e}$ gives back the open 
book $(\Sigma, \phi)=M$. Note that if the fixed 
points $z,z'$ of $\hat{\phi}$ are distinct, then in 
general the sections 
$\{ z\}\times S^1$, $\{ z'\}\times S^1$ may not be isotopic 
(as knots) in 
$\hat{M}$ -- they may even represent different elements of 
$\pi_1(\hat{M})$, hence correspond to different open books. 
The marked point $z$ is essential 
for recovering the original open book and hence contact structure. 
} 
\end{rem} 

\vspace{.2in} 
\noindent
{\bf Extension of $d\alpha$ over $\hat{M}$.} \ 
The 2-form $\hat{\omega}$ canonically extends to a 2-form 
on $\hat{\Sigma}_{\hat{\phi}}$ which we still denote by 
$\hat{\omega}$. 
The differential
$d\alpha|_{\Sigma_\phi}$ extends to a 
closed 2-form $\hat{\tau}$ on  
$\hat{M}$ such that 
$\hat{\tau}$ pulls back to $\hat{\omega}$ on each 
fiber of $\hat{\pi}$. 
Such $\hat{\tau}$ can be expressed as 
\[ 
\hat{\tau} =\hat{\omega}+\eta_t \wedge dt, 
\]   
where $t\in S^1=\R/\Z$ is the coordinate of the base $S^1$,
$\eta=\eta_t$  is a family of closed 1-forms on $\hat{\Sigma}$ 
satisfying $\hat{\phi}^*\eta_0=\eta_1$.

\vspace{.2in} 
\noindent
{\bf Orbits.} \ 
The horizontal distribution $\ker{\hat{\tau}}$ 
is generated by the vector field 
\[ 
\hat{R}:=\hat{\partial_t} + X_\eta,  
\]   
with $\hat{\partial_t}\subset \ker\hat{\omega}$, 
$\hat{\pi}_*(\hat{\partial_t})=\partial_t$, and  
$X_\eta\subset\ker\hat{\pi}_*$ is the $t$-dependent symplectic 
vector field defined by the equation 
\[  
\hat{\omega}(X_\eta, \cdot )=\eta_t.  
\]   

\begin{rem} 
{\rm 
Similar to Remark \ref{Rorbit}, $\hat{R}$-orbits 
correspond to periodic points of some $\hat{\phi}'\in
\text{Symp}(\hat{M},\hat{\omega})$ symplectically isotopic to 
$\hat{\phi}$. 
} 
\end{rem} 

\begin{rem} 
{\rm 
The 2-form $\hat{\tau}$ is the canonical extension of 
$\hat{\omega}\in\Omega^2(\hat{\Sigma})$ over
$\hat{\Sigma}_{\hat{\phi'}}\cong \hat{\Sigma}_{\hat{\phi}}$ 
corresponding to the monodromy $\hat{\phi'}$. 
} 
\end{rem}

\vspace{.2in} 
\noindent 
{\bf Periodic Floer homology.} \ 
Based on the idea of Seidel \cite{Se}, 
Hutchings and Sullivan \cite{HS} 
defined {\em Periodic Floer homology} for $\hat{\phi'}
\in\text{Symp}(\hat{\Sigma},\hat{\omega})$ by counting 
in the symplectic manifold 
\[ 
(\R\times \hat{M},\hat{\tau}+ds\wedge dt), \quad s\in \R, 
\] 
pseudoholomorphic curves converging to periodic trajectories of 
$\hat{R}$ at $s=\pm\infty$. The relevant almost 
complex structure $\hat{J}$ will be 
\begin{enumerate} 
\item $\R$-invariant, 
\item {\em tamed} by $\hat{\Omega}:=\hat{\tau}+ds\wedge dt$,  
i.e. $\hat{\Omega}(v,\hat{J}v)\neq 0$ for $v\neq 0$, and 
\item 
$\hat{J}(\partial_s)=f\hat{R}$, $\hat{J}(f\hat{R})=-\partial_s$ for
some positive function $f$ on $\hat{M}$. 
\end{enumerate}

\vspace{.2in} 
\noindent 
{\bf A motivation.} \ 
Now that each of $(M,\xi)$ and $\hat{M}$ has 
its own theory based on 
counting holomorphic curves, and these two theories have some very 
similar ingredients, namely periodic trajectories and 
almost complex structures. In fact we can have 
\[ 
R=f\hat{R}\quad \text{and}\quad  J=\hat{J} \quad \text{ on } 
\Sigma_\phi. 
\] 

So if we would like to compute (cylindrical) contact homology of 
$(M,\xi)$, and since in doing so we need to know how to count 
holomorphic curves converging to $B^m$ or intersecting with 
$\R\times B$, it helps to know how such curves correspond to
$\hat{J}$-holomorphic curves in $\R\times\hat{M}$, if such a
correspondence does exist.  Moreover, 
One can use the correspondence to
compare the two holomorphic  curve theories, exploring the 
relation between contact homology theory and mapping class groups. 

Our goal in the next section is to establish a holomorphic 
0-surgery model that allows a direct comparison of holomorphic 
curves before and after the surgery.

\section{A holomorphic 0-surgery model} \label{model} 

In this section we construct a holomorphic model of the 
canonical 0-surgery along the binding of an open book. 
It is no news that one can describe a 0-surgery in $\C^2$. 
The novelty here is to do it carefully enough and to find a 
nice vector field $Y$ whose flow (i) preserves the standard 
complex structure $J_o$ on $\C^2$, and (ii) 
embeds the two symplectizations into $\C^2$ nicely so that 
$J_o$ is their common almost complex structure (see Lemma 
\ref{symps}).   This is done in Sections \ref{modelN}-\ref{2symp}. 
Lemma \ref{Nview} and \ref{N'view} 
in Section \ref{holo} 
describe how punctured holomorphic discs in $\C^2$ are 
perceived in each symplectizations. Corollary \ref{inter} 
in particular will be used to determine 
the boundary operator (see Section \ref{holocyl}). 
In Section \ref{global} it is confirmed that correspondences  
between holomorphic curves in the local model extend straight 
forwardly to correspondence between pseudoholomorphic 
curves in $\R\times M$ and $\R\times \hat{M}$.

\subsection{A contact solid torus in $\C^2$} \label{modelN}

Let $\C^2:=\{ (z_1,z_2)\mid z_1,z_2\in\C\}$ be the complex plane 
with the standard complex structure $J_o$. Write 
$z_j=r_je^{i\theta_j}$. 

Consider on $\C^2$ the smooth function  
\[  
F(z_1,z_2):=-|z_1|^2+|z_2|^2=-r_1^2+r_2^2. 
\]  
Let $\epsilon >0$ be a constant and define 
\begin{equation} \label{N} 
N=N_{\epsilon}:=F^{-1}(-1)\cap \{ r_2<\epsilon\}. 
\end{equation}  
$N$ is diffeomorphic to $S^1\times D^2$. 
Let 
$\iota:N\hookrightarrow \C^2$ denote the 
inclusion map. Define 
\[ 
\lambda:=\iota^*(\frac{-1}{2}dF\circ J_o )
=-(1+r_2^2)d\theta_1+r_2^2d\theta_2.  
\] 
\begin{fact} 
The 1-form $\lambda$ is a contact on $N$, its contact structure 
and  Reeb vector field are 
\[ 
\zeta:=\ker\lambda=TN\cap J_oTN, \qquad 
R_\lambda=-\partial_{\theta_1}-\partial_{\theta_2}
\] 
\end{fact} 
Then $\zeta$ is the maximal complex subbundle of the tangent bundle of 
$N$. $(N,\zeta)$ will be a model of a tubular neighborhood of the
binding.  Note that all trajectories of $R_\lambda$ are periodic. 
We would like to perturb $\lambda$, but keeping $\zeta$ intact, 
so that the resulting Reeb vector field has only one simple 
periodic orbit namely, $\gamma:=\{ r_1=1\} \cap N$ with orientation 
given by $-\partial_{\theta_1}$.  

\begin{lem} \label{lambda'} 
For any constant $c>0$  
there is a smooth function $h=h(r_2)\in C^\infty(N)$ depending only 
on $r_2$ and $c$ such the the Reeb vector field 
$R_{\lambda'}$ of the contact 1-form 
$\lambda':=e^{-h}\lambda$ is
$R_{\lambda'}=-\partial_{\theta_1}+c\partial_{\theta_2}$. 
In particular, if $c\in\R_+\setminus \Q$ then the only simple 
Reeb orbit of $R_{\lambda'}$ is 
$\gamma:=\{ r_1=1\} \cap N$ with orientation 
given by $-\partial_{\theta_1}$.
\end{lem} 

\begin{proof} 
Let $h=h(r_2)\in C^\infty(N)$ then the Reeb vector field 
of $\lambda':=e^{-h}\lambda$ is 
\[ 
R_{\lambda'}=e^h(R_{\lambda}+Z_h),  
\] 
where $Z_h\subset \zeta$ is the unique vector field satisfying 
\[ 
d\lambda (Z_h,\cdot )=-dh  \quad \text{ on }\zeta.  
\] 
A straightforward calculation yields 
\[  \hspace{.7in} 
Z_h=\frac{h'}{2r_2}(r_2^2\partial_{\theta_1}+
(1+r_2^2)\partial_{\theta_2}), \hspace{.8in}  (h':=\frac{dh}{dr_2}).  
\] 
Then 
\[ 
R_{\lambda'}=e^h\Big( \big(\frac{r_2h'}{2}-1\big)\partial_{\theta_1}
+\big(\frac{(1+r^2_2)h'}{2r_2}-1\big)\partial_{\theta_2}\Big).  
\] 

Fix a positive constant $c$ and solve for 
$h$ satisfying the initial value problem  
\[ 
\frac{\frac{(1+r^2_2)h'}{2r_2}-1}{\frac{r_2h'}{2}-1}=-c, \quad 
h(0)=0  
\] 
We get 
\[ 
h(r_2)=\ln(1+(1+c)r_2^2),  
\] 
and 
\begin{equation} \label{R'} 
R_{\lambda'}=\iota^*(X), \quad 
X:=-\partial_{\theta_1}+c\partial_{\theta_2}. 
\end{equation}  
\end{proof} 

Let 
\[ 
\pi_N:N\setminus \gamma\to S^1_{\theta_2}
\] 
denote the projection 
onto the $\theta_2$-coordinate. Then $(\pi_N,\gamma)$ is an open 
book representation of $N$ with pages diffeomorphic to an annulus. 
The monodromy $\psi$, which is the time 1 map of the flow of 
$\frac{1}{c}R_{\lambda'}$, is isotopic to the identity map.

\begin{cor} 
The 2-form $d\lambda'$ is $\theta_2$-independent and is 
symplectic when restricted to any 
page of $\pi_N$. Let $\omega$ denote the restriction of 
$d\lambda'$ to a page. Then $\psi^*\omega=\omega$. 
\end{cor} 

The following lemma shows that the binding $\gamma$ and its 
positive iterates $\gamma^m$ are {\em elliptic} (see (\ref{z2ind})). 
This result will be used in Section \ref{reeb}. 

\begin{lem} \label{Beven} 
Let $R_{\lambda'}$ be as in (\ref{R'}) with $0<c\not\in \Q$. 
Let $\gamma:=N\cap \{ r_1=1\}$ be the unique simple Reeb orbit on 
$N$. Then for $m\in\N$, 
$\gamma^m$ is elliptic, i.e., $\bar{\mu}(\gamma^m,\Z_2)=0$. 
\end{lem} 

\begin{proof} 
Note that 
\[ 
R_{\lambda'}=x_1dy_1-y_1dx_1+c(x_2dy_2-y_2dx_2). 
\] 
Let $\zeta:=\ker\lambda'$. Since $\zeta|_\gamma=\text{span}
(\partial_{x_2},\partial_{y_2})$, then with respect to the 
ordered basis $\{ \partial_{x_2},\partial_{y_2}\}$, 
\[ 
DR_{\lambda'}|_\gamma=\begin{bmatrix}0&-2c\\ 2c&0\end{bmatrix}.  
\]  
Since the action of $\gamma^m$ is $2m\pi$, the linearized 
Poincar\'{e} return map of the flow of $R_{\lambda'}$ along 
$\gamma^m$ is 
\[ 
\Lambda_{\gamma^m}= 
\exp \Big( 2m\pi \begin{bmatrix}0&-2c\\ 2c&0\end{bmatrix}\Big) 
=\begin{bmatrix}\cos 4mc\pi & -\sin 4mc\pi \\ \sin 4mc\pi 
& \cos 4mc\pi  \end{bmatrix}. 
\] 
If $c\not\in \Q$ then $\det(\Lambda_{\gamma^m}-Id)<0$, so 
$\gamma^m$ is elliptic. 
\end{proof}

\subsection{A 0-surgery in $\C^2$} \label{0-surg}

Recall the small constant $\epsilon >0$ and 
\[ 
N=N_{\epsilon}:=F^{-1}(-1)\cap\{ r_2<\epsilon\}.  
\] 

%
%

%
%

%
%

On $N$ we apply the canonical 0-surgery along $\gamma$ to get a 
new manifold 
\[ 
\hat{N}:=\{ -e^{-2s}r_1^2+e^{2cs}r_2^2=1\}\cap \{ r_2<\epsilon\} 
\] 
satisfying the following conditions:  

\begin{enumerate} 
\item $s=s(r_2)>0$ is a smooth function depending only on $r_2$, 
\item $\hat{N}\pitchfork Y$ where $Y:=-J_oX=-r_1\partial_{r_1}
+cr_2\partial_{r_2}$, 
\item $\hat{N}\cap \{ r_2>\epsilon'\} =N\cap\{ r_2>\epsilon'\} $ 
for some positive constant $\epsilon'<\epsilon$.  
\end{enumerate} 

Note that $\hat{N}$ is a solid torus diffeomorphic to 
$D^2\times S^1_{\theta_2}$. 
Let 
\[ 
\pi_{\hat{N}}:\hat{N}\to S^1_{\theta_2}
\] 
denote the projection onto the $\theta_2$-coordinate. 
$\pi_{\hat{N}}$ gives $\hat{N}$ the structure of a disc bundle over 
$S^1_{\theta_2}$. 

The vector field $X=-\partial_{\theta_1}+
c\partial_{\theta_2}$ is tangent to $\hat{N}$ 
and transversal to the fibers of $\pi_{\hat{N}}$. 
Together $\pi_{\hat{N}}$ and 
$X$ induce on $\hat{N}$ the structure of a mapping torus 
of a disc with monodromy $\hat{\psi}$ induced by the time 1 map 
of the flow of $\frac{1}{c}X$.

Let $\Omega:=r_1dr_1\wedge d\theta_1+r_2dr_2\wedge d\theta_2$. 
$\Omega$ is the standard symplectic 2-form on $\C^2$. 
$\Omega$ is invariant under the flow of $X$ and pulls back to 
a symplectic 2-form on every page of $\pi_{\hat{N}}$ (as well as 
on every fiber of $\pi_N:N\setminus \gamma\to S^1_{\theta_2}$). 
$\hat{\psi}$ is symplectic with respect to this pulled back 2-form, 
and is isotopic to the identity map. $\hat{\psi}$ has only one 
fixed point, which corresponds to the simple loop 
\[ 
\hat{\gamma}:=\hat{N}\cap \{ r_1=0\}, \quad \dot{\gamma}=\frac{1}{c}X.  
\]

\subsection{Two overlapping symplectizations in $\C^2$}
\label{2symp}

Recall the vector field 
\[ 
Y:=-J_oX=-r_1\partial_{r_1}+cr_2\partial_{r_2}. 
\]

Let $Y^t$, $t\in \R$, denote the flow of $Y$. 
\begin{fact} 
The flow $Y^t$ preserves $J_o$ and $X$. 
\end{fact}

Note that $Y^t$ also preserves the values 
$\rho:=r_1^cr_2$, $\theta_1$ and $\theta_2$. 
The integral trajectories of $Y^t$ on $\C^2\setminus\{ 0\}$ 
are parametrized by $(\rho,
\theta_1,\theta_2)$, $\rho>0$, $\theta_1,\theta_2\in S^1=\R/2\pi\Z$.

Recall from Lemma \ref{lambda'} the contact 1-form $\lambda'$ on 
$N$. The contact structure $\zeta=\ker\lambda'$ is a $J_o$-complex 
line bundle, $J_o|_\zeta$ is $d\lambda'$-compatible, 
and $J_oY|_N=X|_N$ is the Reeb vector field $R_{\lambda'}$. 

With the above understood we can embed the symplectization 
$\R\times N$ into $\C^2$ by 
\begin{equation} 
\text{identifying } \{ t\}\times N \text{ with } 
Y^tN , \quad \forall t\in  \R,  
\end{equation}  
where $Y^tN$ is the image of $N$ under 
the time $t$ map of the flow of $Y$. In particular the vector 
field $\partial _t$ of $\R\times N_\epsilon$ is identified with 
$Y$. 
Let 
\begin{equation}  
W:=\cup_{t\in\R}Y^tN \cong \R\times N. 
\end{equation}  
Observe that $J_o$ is $Y$-invariant hence is a 
$\lambda'$-admissible almost 
complex structure on $\R\times N\cong W$. 

Likewise, we can also embed the "symplectization"  
$\R\times \hat{N}$ into $\C^2$ via the flow of $Y$ 
 by 
\begin{equation} 
\text{identifying } \{ t\}\times \hat{N}\text{ with } 
Y^t\hat{N} , \quad \forall t\in  \R.  
\end{equation}  
Now the vector field $\partial_t$ of $\R\times\hat{N}$ is 
identified with $Y$ as well. Let 
\begin{equation}  
\hat{W}:=\cup_{t\in\R}Y^t\hat{N} 
\cong\R\times\hat{N}.  
\end{equation}  
Again, $J_o$ is $Y$-invariant and $\Omega$-compatible.

\begin{lem} \label{symps}
The symplectizations $\R\times N$ and $\R\times 
\hat{N}$ are overlapping on the region 
$W\cap\hat{W}\subset \C^2\setminus \{ z_1z_2=0\}$, 
and they share the same compatible complex structure $J_o$. 
Hence a $J_o$-holomorphic map into $W\cap\hat{W}$ 
is a $J_o$-holomorphic map into both symplectizations. 
\end{lem}

\subsection{Holomorphic maps into $\R\times N$ and 
$\R\times \hat{N}$} \label{holo}

Recall that the binding of $N$ is the 
Reeb orbit 
\[ 
\gamma:=N\cap \{r_2=0\} \text{ oriented by } X|_{\gamma}.    
\] 
The 0-surgery along $\gamma$ has the 
effect of replacing the open book $N$ by the 
mapping torus $\hat{N}$, and the binding 
$\gamma$ by the orbit 
\[ 
\hat{\gamma}:=\hat{N}\cap \{|z_1|=0\} \text{ oriented by }
(X)|_{\hat{\gamma}}. 
\] 

With the notations 
\[ 
\C^*_{z_1}:=(\C\setminus\{ 0\})\times \{ 0\}, \qquad 
\C^*_{z_2}:=\{ 0\}\times (\C\setminus\{ 0\} ), 
\] 
we have 
\begin{gather} 
\R\times \gamma\cong \cup_{t\in \R}Y^t\gamma=
 \C^*_{z_1},\\ 
\R\times \hat{\gamma}\cong \cup_{t\in \R}Y^t\hat{\gamma}=
\C^*_{z_2}. 
\end{gather} 
The union of both is then contained in the set 
\[ 
z_1z_2=0.  
\]

Let $U\subset \C$ be an open disc containing the point $z=0$. 
Denote $U^*:=U\setminus\{ 0\}$. We are interested in $J_o$-holomorphic 
maps 
\[ 
f(z)=(f_1(z),f_2(z)):=U^*\to W\cap\hat{W}
\subset\C^2 \setminus \{ z_1z_2=0\} 
\] 
satisfying 
\begin{cond} \label{condi} 
$f(z)$ converges to a positive multiple of either 
$\gamma$ or $\hat{\gamma}$ at either $t=\infty$ or $t=-\infty$ 
as $z\to 0$.  
\end{cond} 

Recall that $Y^t$ preserves the value $r_1^cr_2$, so two points 
$(z_1,z_2), (z'_1,z'_2)$ 
in $\C^2$ are projected by $Y^t$ to the same point either in 
$N$ or in $\hat{N}$ only if $|z_1|^c|z_2|=|z'_1|^c|z'_2|$. 
Since $\R\times\gamma\cup \R\times\hat{\gamma}$ is contained 
in the set $\{ z_1z_2=0\}=\{r_1^cr_2=0\}$, then a sequence of points 
$z_k=(z_{1k},z_{2k})$ in $N$ (resp. in $\hat{N}$) approach  
$\gamma$ (resp. $\hat{\gamma}$) implies that 
$z_{1k}^cz_{2k}\to 0$ as $k\to \infty$. 

Condition \ref{condi} 
then implies that 
\begin{equation}  
f_1(z)^cf_2(z)\to 0 \quad \text{ as }z\to 0.  
\end{equation} 
Hence for $j=1,2$, either $f_j$ can be 
holomorphically extended over $z=0$ or $f_j$ has a pole of 
finite order at $z=0$. Since 
\[ 
 Y^t(z_1,z_2)=(e^{-t}z_1,e^{ct}z_2) 
\]
so if $|f_j|\in O(|z|^{n_j})$ near $z=0$ 
for some $n_j\in \Z$, $j=1,2$, then 

\begin{equation}  \label{cn1n2} 
cn_1+n_2>0.  
\end{equation} 

\begin{rem} 
Often we will take $c\in \R_+\setminus\Q$ to be arbitrarily large, 
then (\ref{cn1n2}) eventually implies that $n_1\geq 0$ 
(when $c\to \infty$). 
\end{rem}

The sign of $n_1$ (resp. $n_2$) 
determines how $f$ is perceived from the point of 
view of $N$ (resp. $\hat{N}$). 

\begin{lem} \label{Nview} 
View $f$ as a $J_o$-holomorphic map into the symplectization 
$\R\times N$. Then we have the following conclusions  
depending on the sign of $n_1$. 

\begin{enumerate} 
\item $n_1>0$. $f(z)$ converges to $\gamma^{n_1}$ at
$t=\infty$  asymptotically as $z\to 0$. 

\item $n_1=0$. $f$ can be holomorphically extended over $z=0$. 
$f(U)$ intersects transversally and positively with 
$\C^*_{z_1}$. 
The intersection multiplicity is $n_2$.

\item $n_1<0$. $f(z)$ converges to $\gamma^{-n_1}$ at
$t=-\infty$  asymptotically as $z\to 0$.
\end{enumerate} 
%
%
\end{lem} 

By interchanging  $n_1$ and $n_2$ we get similar conclusion from the 
perspective of $\hat{N}$. 

\begin{lem} \label{N'view} 
View $f$ as a $J_o$-holomorphic map into the symplectization 
$\R\times \hat{N}$. Then we have the following conclusions 
depending on the sign of $n_2$. 

\begin{enumerate} 
\item $n_2>0$. $f(z)$ converges to 
$\hat{\gamma}^{n_2}$ at
$t=-\infty$  asymptotically as $z\to 0$. 

\item $n_2=0$. $f$ can be holomorphically extended over $z=0$. 
$f(U)$ intersects transversally and positively with 
$\C^*_{z_2}$. 
The intersection multiplicity is $n_1$.

\item $n_2<0$. $f(z)$ converges to $\hat{\gamma}^{-n_2}$ at
$t=\infty$ asymptotically as $z\to 0$.
\end{enumerate} 
%
%
\end{lem}

When computing contact homology of an open book, 
one needs to take into account holomorphic curves intersecting 
with the holomorphic cylinder $\R\times B$, where $B$ is the 
binding of the open book. In our local model here, $\R\times B$ is 
identified with $R\times \gamma=\C^*_{z_1}$. The following 
corollary, as a special case of Lemma \ref{mainlem}, states that, in our local model, 
how such curves correspond to curves in $\R\times \hat{N}$. 

\begin{cor} \label{inter} 
There is a one-one 
correspondence between pseudoholomorphic discs in $\R\times N$ 
that intersect with $\R\times\gamma$ at one point 
with winding number $n_2>0$ and 
pseudoholomorphic half-cylinders in $\R\times \hat{N}$ that 
converges to $\hat{\gamma}^{n_2}$ at $t=-\infty$. 
\end{cor}

\subsection{From local to global} \label{global}

Recall the mapping torus $M=\Sigma_\phi\cup_{id}B\times D^2$, 
where 
$B\times D^2$ is  diffeomorphic to a solid torus $S^1\times D^2$. 
Let $p$ be the angular coordinate of $B\cong \R/2\pi\Z$ and 
$(r,\theta)$ the polar coordinates of $D^2=\{ r<\epsilon\}$.

Recall $N\subset \C^2$ from (\ref{N}) as well as the 
contact 1-form $\lambda'$ and its Reeb vector field 
$R_{\lambda'}$ from Section \ref{modelN}. 
The diffeomorphism $\Phi:B\times D^2\to N$ defined by
\begin{equation} \label{Phi} 
\Phi (p,r,\theta)=(\sqrt{1+r^2}e^{-ip},re^{i\theta}) 
\end{equation}  
induces on $B\times D^2$ the contact 1-form $\Phi^*\lambda'$ 
which can be extended over $\Sigma_\phi$ to be a contact 1-form  
$\alpha$ on $M$ whose contact structure $\xi$ is supported by 
the open book $(\Sigma,\phi)$. 

On $B\times D^2$ we have 
$R_\alpha=\Phi^*(X|_N)$ hence 
$B^m=\Phi^{-1}(\gamma^m)$, $m\in\N$, are 
the only orbits of $R_\alpha$ on $B\times D^2$. 

With $B\times D^2$ identified with $N$, the mapping torus 
$\hat{M}$ obtained by a 0-surgery along $B$ can be 
identified with 
\[ 
\hat{M}=\Sigma_\phi\cup_{\Phi|_{B\times S^1}}\hat{N} =
\hat{\Sigma}_{\hat{\phi}}
\] 
where $\hat{\Sigma}=\Sigma\cup D^2$ is the closed Riemann 
surface obtained by gluing along $\partial \Sigma\cong S^1$ 
a 2-disc with 
$\hat{\phi}|_{\Sigma}=\phi$ and $\hat{\phi}=id$ on 
$\hat{\Sigma}\setminus \Sigma$. 
The Reeb vector field $R_\alpha$ induces a vector field 
$\hat{R}$ on $\hat{M}$ such that 
$\hat{R}|_{\Sigma_\phi}=R_\alpha$ and $\hat{R}|_{\hat{N}}
=X|_{\hat{N}}$.

With the diffeomorphism $\Phi$, 
the common almost complex structure $J_o$ of $\R\times N$ and 
$\R\times \hat{N}$ extends over $\Sigma_\phi$ as an 
$\alpha$-admissible almost complex structure. 
The holomorphic curve correspondence between $\R\times N$ 
and $\R\times \hat{N}$ now extends over $\R\times M$ and 
$\R\times \hat{M}$ straightforwardly. 
In particular, 
Lemma \ref{Nview} and Lemma \ref{N'view} together imply Lemma  
\ref{mainlem}.

\section{Cylindrical contact homology of a Dehn twist} 
\label{appl}

We define the Dehn twist in Section \ref{Dehn} and show in 
Section \ref{reeb} that, on the 
corresponding open book $M$, the Reeb orbits come in three different 
types (Proposition \ref{bct}). Holomorphic cylinders are studied 
in Section \ref{holocyl}. It is shown that under various topological 
and dimensional constraints, the cylindrical contact homology is 
defined for $M$, and 
there are only a handful of types of 
holomorphic cylinders to be counted. Lemma \ref{partial} gives a 
complete classification (up to sign) of 1-dimensional moduli of 
holomorphic cylinders, either intersecting with the binding or not. 
An energy estimate is calculated in Section \ref{energy-est} 
to help understand the coherent orientation of the moduli 
and hence the signs involved in the boundary operator of the 
cylindrical contact homology $HC(M,\xi)$ 
(Section \ref{orien}). We complete the computation of $HC(M,\xi)$ 
in Section \ref{hc}, the result  is summarized in Theorem \ref{example}.

\subsection{Dehn twist} \label{Dehn} 

Let $\Sigma:=T^2\setminus D^2$ be a 2-torus with a disc removed. 
Let $\Gamma\subset \Sigma$ be an embedded nonseparating circle, 
$U\subset\overset{\circ}{ U'}, U'\subset\Sigma$ be 
closed tubular neighborhoods of $\Gamma$. Let $(q,p)$ be the 
coordinates  of $U'\cong [q'_-,q'_+]\times S^1$, 
$U\cong [q_-,q_+]\times S^1$, $q'_-<q_-<q_+<q'_+$,  
$S^1\cong\R/\Z$, 
$\Gamma\cong\{ q_o\}\times S^1$ for some $q_o\in (q_-,q_+)$.  
We also use coordinates $(q,p)$ near $\partial\Sigma$ so that 
$\partial\Sigma=\{ q=const\}$ and $dq\wedge dp$ is an area form 
near $\partial\Sigma$. We fix a symplectic 2-from $\omega$ 
on $\Sigma$ so that 
\[ 
\hspace{1in} 
\omega =d(qdp)  \quad \text { on $U'$ and near $\partial\Sigma$}. 
\] 

Fix a natural number $\sigma\in\N$ and let $\phi\in\text{Symp}
(\Sigma,\omega )$ denote a $\sigma$-Dehn twist 
supported on $U\supset \Gamma$ such that 
\begin{align*} 
\phi &=id \hspace{.83in} \text{ on } \Sigma\setminus U, \\ 
\phi (q,p)&=(q,p-f(q)) \quad \text{ on } U=[q^-,q^+]\times S^1, 
\end{align*} 
with $f:[q^-,q^+]\to [0,\sigma ]$ a smooth surjective 
increasing function satisfying 
\[ 
 f'(q) \geq 0 , \quad f'(q^{\pm})=0.  
\]

The 
mapping torus $U_\phi$  can be described topologically as 
follows: Let $T^2=\R /\Z \times\R /\Z$ be a torus parametrized by 
$(p,t)$ and let $[p]$, $[t]$ be the corresponding generators of 
$\pi_1(T^2)$. 
Take two thickened tori $[q_-,q_o]\times T^2$ and $[q_o,q_+]\times T^2$ 
and glue them together by using the map $g_\tau: \{ q_o\}\times T^2\to 
\{ q_o\}\times T^2$, 
\begin{equation*} 
g_\tau (p,t)= (p,t+\sigma p).  
\end{equation*} 

\begin{prop} \label{-htpy+} 
Write $\partial U_\phi=T_-\cup T_+$, where 
$T_-:=\{ q_-\}\times S^1_p\times S^1_t$ 
and $T_+:=\{ q_+\}\times  S^1_p\times S^1_t$. Let 
$[p]=[\Gamma]$ and $[t]$ denote the generators of 
$T_{\pm}\cong S^1_p\times S^1_t$ 
with respect to the coordinates $(p,t)$. 
Then as elements of $\pi_1(U_\phi)$ we have 
\begin{equation}  \label{htpy} 
\pi_1(T_-)\ni [p]^n[t]^m=
[p]^{n+\sigma m}[t]^m\in \pi_1(T_+).  
\end{equation} 
\end{prop} 

\begin{proof} 
Consider the homotopy $H:[q_-,q_+]\times S^1_\tau\to [q_-,q_+]
\times S^1_p\times S^1_t$,  
\[ 
H(q,\tau):=(q,n\tau +m\tau f(q) , m\tau ).  
\] 
\end{proof}

\subsection{Contact 1-form and Reeb orbits} \label{reeb}

Let $M$ denote the 3-manifold represented by the open book 
$(\Sigma,\phi)$ as defined in Section \ref{Dehn} 
and let $B$ denote the binding. In this section we will 
construct an explicit contact 1-form $\alpha$ associated to 
the open book $(\Sigma,\phi)$ and classify its Reeb orbits.

\vspace{.2in} 
\noindent
{\bf Constructing $\alpha$.} \ 
Let $\beta\in\Omega^1(\Sigma)$ be a primitive of $\omega$ 
such that $\beta=qdp$ on $U$ and near $\partial\Sigma$. 
We have $\phi^*\beta-\beta=-qf'(q)dq$. Define 
\[ 
\alpha:=\beta -tqf'(q)dq +Kdt,  
\] 
with 
$K\in C^\infty (\Sigma)$ to be specified in Condition \ref{K} below. 
It is easy to see that  $\alpha$ 
is contact  if $K$ is a positive constant large enough. 
Note that $\alpha =\beta +Kdt$ away from the Dehn twisted regions. 

If $K$ is a constant then the Reeb vector field $R$ on the mapping torus 
$\Sigma_\phi$ is $R=(\partial _t - qf'\partial _p)/(K -q^2f' )$. 
On $U \times [0,1]$ 
the flow of $R$ takes a point $(q,p)$ at $t=0$ to the point 
$(q,p-qf'(q))$ at $t=1$, then $\phi$ identifies 
$(q,p-qf'(q))$ at $t=1$ with 
the point 
\[ (q,p-qf'(q)-f(q))=(q,p-(qf(q))') \quad \mbox{ at } t=0. 
\] 
The coordinate $q$ is fixed under the flow of $R$. 
When $(qf(q))' \in \Q$  the $q$-level 
set, which is a 2-torus, is fibrated by $R$-orbits; while 
there are no periodic Reeb trajectories if $(qf(q))'\not\in\Q$. 
So we get infinitely many $S^1$-families of $R$-orbits on 
the mapping torus $U _\phi$. 
Moreover $(\Sigma \setminus U )_{\phi}=
(\Sigma\setminus U)\times S^1_t$  
is  a trivial $S^1$-bundle with Reeb orbits as fibers. Such 
a contact 1-form $\alpha$ is not regular. Recall that $\alpha$ is 
regular if the set of Reeb orbits is discrete, otherwise $\alpha$ is 
called not regular. 

Instead of letting $K$ be 
a constant we want to choose 
$K$ so that the resulting Reeb vector field has no 
periodic trajectories in $U_\phi$. Let $K=K(q)$ on $U$, then $K$ is 
$\phi$-invariant hence can be thought as a function on $U_\phi$ that 
depends only on $q$. The Reeb vector field is then 
$(\partial_t-(qf'+K_q)\partial_p)/(K-q(qf'+K_q))$. 
Level $q=const$ has 
$R$-orbits if and only if $qf'+K_q+f=K_q+(qf)'\in \Q$. So if 
$K_q=-(qf)'-\tilde{r}$ for some $\tilde{r}\not\in\Q$ then there 
will be no $R$-orbits on $U_\phi$. Such $K$ is equal to 
$-q(f+\tilde{r})+c_0$ for some constant $c_0$. To ensure  $\alpha$ 
being contact we need $K-qK_q-q^2f>0$, which implies that $c_0>0$. 

Now consider a smooth function $K\in C^{\infty}(\Sigma)$ 
satisfying the following  

\begin{cond} \label{K} 
\begin{enumerate} 
\item On $U$, $K=K(q)=-q(f+\tilde{r})+c_0$ for some constant $c_0>0$, 
$\tilde{r}\not\in \Q$, $\tilde{r}\lesssim 0$, i.e., $\tilde{r}<0$ and 
$\tilde{r}\sim 0$. 

\item $K=K(q)$ on $U'$; $K_q>0$ for $q\in (q'_-,q_-)$; 
$K_q<0$ for $q\in (q_+,q'_+)$; and $K_{qq}>0$ on $U'\setminus U$. 

\item $K$ is a Morse function on $\Sigma\setminus U$, with 
only one critical point $x_h$ which is hyperbolic; $|dK|\sim 0$ on 
$\Sigma\setminus U'$.  

\end{enumerate} 
\end{cond} 
Moreover, we can extend $K$ over $B\times D^2$ so that on $B\times D^2$ 
(see (\ref{Phi}) and Lemma \ref{lambda'}) 
\[ 
\alpha =q(r)dp+K(r)d\theta =\Phi^*\lambda ',  
\] 
which will ensure that the only Reeb orbits of $\alpha$ on $B\times
D^2$ are $B^m$, $m\in\N$. 

Let $R=R_\alpha$ denote the corresponding Reeb vector field, 
then 
\begin{equation} \label{R-phi} 
R=
\begin{cases} 
(\partial_t+(f+\tilde{r})\partial_p)/(K-qK_q-q^2f)& 
\text{ on } U_\phi,\\  (\partial_t-K_q\partial_p)/(K-qK_q)& \text{ on }
(U'\setminus U)_\phi, \\  
(\partial_t-X_K)/(K-\beta (X_K))& \text{ on } (\Sigma\setminus 
U')_\phi,  
\end{cases} 
\end{equation} 
where $X_K$ is the Hamiltonian vector field of $K$ relative to the 
symplectic 2-form $\omega =d\beta$, 
i. e., $d\beta (X_K,\cdot )=-dK$.

\vspace{.2in} 
\noindent
{\bf Types of Reeb orbits.} \ 
The following proposition describes the three types (
Type $\fT$, $\fC$ and $\fB$) 
of  Reeb orbits on $(M,\alpha)$: 

\begin{prop} \label{bct} 
There are three types of Reeb orbits described as follows: 
\begin{enumerate} 
\item Type $\fT$ orbits. These orbits come from the 
$\sigma$-Dehn twist and lie in $U'\setminus U$. 
They are parametrized by the set 
\begin{equation} \label{nm} 
\sS:=\{ (n,m)\in \N\times\N \mid   0<n<\sigma m\}. 
\end{equation} 
(See Remark \ref{aboutT} below for more discussion on the 
parametrization.) 
When $K$ unperturbed, each of such pair 
$(n,m)$ represents an $S^1$-family of 
(degenerate) Reeb orbits contained in the 2-torus 
$\{ K_q=\frac{-n}{m}\} \cap U'$. 

A further perturbation of $K$ will turn each $S^1$-family of orbits 
into a pair of orbits denoted by $h_{n/m}$ and $e_{n/m}$, 
$h_{n/m}$ is hyperbolic while $e_{n/m}$ is elliptic (see \cite{HS}). 

\item Type $\fC$ orbits. They are 
$h^m$, $m\in \N$, where $h=\{ x_h\}\times S^1_t$ 
is the simple Reeb orbit corresponding to the unique critical 
point $x_h$ of $K|_{\Sigma\setminus U}$, and $x_h$ is hyperbolic.

\item Type $\fB$ orbits. They are $B^m$, $m\in \N$.  
\end{enumerate} 
\end{prop} 

\begin{rem} \label{aboutT}
{\rm 
Recall the number $\tilde{r}\lesssim 0$, 
$\tilde{r}\not\in\Q$, from Condition \ref{K}. 
The $\fT$-orbits in the region $(q_+,q'_+)$ are parametrized by the 
set 
\[ 
\sS_+:=\{ (n,m)\in \N\times\N\mid  \ n<(\sigma+\tilde{r})m\}, 
\] 
where $(n,m)\in H_1(T_+=S^1_p\times S^1_t, \Z)$ is indeed the 
homology class of the corresponding orbits. Similarly  
the $\fT$-orbits in the region $(q'_-,q_-)$ are parametrized 
by 
\[ 
\sS_-:=\{ (n,m)\in (-\N)\times\N\mid  \tilde{r}m<n<0\},  
\] 
where $(n,m)\in H_1(T_-=S^1_p\times S^1_t, \Z)$.  

Recall from (\ref{-htpy+}) that an element $(n,m)\in H_1(T_-,\Z)$ 
gets identified with $(n+\sigma m,m)\in H_1(T_+,\Z)$. Thus $\sS_-$ 
is identified with 
\[ 
\{ (n,m)\in\N\times\N\mid  (\sigma+\tilde{r})m<n<\sigma m\}, 
\] 
where $(n,m)\in H_1(T_+=S^1_p\times S^1_t, \Z)$. 
As a result, 
$\fT$-orbits are parametrized by their corresponding homology classes 
in $T_+$, i.e., by the set $\sS$ in (\ref{nm}).  
 
In fact, as $\tilde{r}$ can be arbitrarily close to 0, when 
$\tilde{r}\to 0$, all $\fT$-orbits are "pushed into" 
the region $q_+<q<q'_+$ and hence are parametrized by $\sS$.  
} 
\end{rem} 

\begin{rem} 
{\rm 
Indeed there are Reeb orbits other than the three types discussed 
above. These extra orbits 
lie in the region between $x_h$ and the collar of $B$, 
and all wind around $B$ with winding numbers $>c$ (recall the number 
$c$ from Lemma \ref{lambda'}), and are homotopic to $B^m$ for some 
$m\in\N$. Since $c$ can 
be made arbitrarily large, these "transient" orbits disappear 
as $c\to\infty$. Moreover, as these orbits are homologously trivial 
and $c_1(\xi)=0$ (see Lemma \ref{c1=0} below), 
their $\bar{\mu}$-indexes are defined and will become arbitrarily
large as $c\to\infty$, 
hence they have no contribution to the 
contact homology. The only Reeb orbits responsible for the 
contact homology are of Type $\fT$, $\fC$ and $\fB$. 
} 
\end{rem}

\begin{lem} \label{allgood} 
All Reeb orbits are {\em good}. 
\end{lem} 

\begin{proof} 

First of all, from Lemma \ref{Beven} one sees that for all $m\in\N$,
$B^m$ is elliptic and its Poincar\'{e} return map has no real
eigenvalues (assuming $c\not\in\Q$). So 
$B^m$ is always good. 

Secondly, $h^m$ is hyperbolic for all 
$m\in \N$, and both eigenvalues of its Poincar\'{e} return map 
are positive, so $h^m$ is good for all $m\in\N$. 

Finally, 
prior to a further perturbation of $K|_{(U'\setminus U)_\phi}$ we have 
$S^1$-families  of Type $\fT$-orbits indexed by $(n,m)$, and 
on $(U'\setminus U)_\phi$ 
\begin{equation} \label{DR} 
DR 
=\frac{K}{(K-qK_q)^2}\begin{bmatrix}0&0\\-K_{qq}&0\end{bmatrix} 
=\begin{bmatrix}0&0\\-\epsilon &0\end{bmatrix}, \quad 0<\epsilon\ll 1.  
\end{equation}  
A slight perturbation of $K$ will deform each $S^1$-family of 
$\gamma_{n/m}$ into a pair of Reeb orbits: $e_{n/m}$ and $h_{n/m}$. 
One can show that the Poincar'{e} return maps of $e_{n/m}$'s 
are negative rotations, hence 
has no real eigenvalues, and the Poincar'{e} return maps for 
$h_{n/m}$'s are hyperbolic with two positive eigenvalues. 
So $e_{n/m}$ and $h_{n/m}$ are good as well. 
\end{proof}

\vspace{.2in} 
\noindent
{\bf Chern class of $\xi$ and contractible orbits.} \ 
Identify $\Sigma$ with the page zero $\Sigma \times \{ 0\}$ of the 
mapping torus $\Sigma _\phi$. 
Let $z\in\Sigma\setminus U$ be a fixed point of $\phi$ 
with  $\{ z\}\times S^1_t$  contractible in $M$. Let
$t_z:=[\{ z\}\times S^1_t]\in\pi_1(\Sigma_\phi,z)$.    
Let $\phi_\#$  be the map on $\pi_1
(\Sigma,x)$ induced by $\phi$. We have the following proposition 
concerning the topology of $\Sigma_\phi$ and $M$. 

\begin{prop} \label{pi1} 
\begin{itemize} 
\item $\pi_1(\Sigma_\phi )=\langle x\in\pi_1(\Sigma,z), t_z\mid
\phi_\#x= t_z^{-1}xt_z \rangle$ 

\item 
$\pi_1(M) =\langle x\in \pi_1(\Sigma,z)\mid \phi_\#x=x\rangle  
= \langle a, b \mid b^\sigma =1\rangle$

\item $H_1(M,\Z)\cong \Z\oplus \Z_{\sigma}$. 

\end{itemize} 
\end{prop}

Note that $[\Gamma]=b$ and 
$[\{ z\}\times S^1_t]=t_z$ 
generates an abelian subgroup $\Z\times \Z$ of $\pi_1(\Sigma_\phi)$, 
as well as an $\Z_{\sigma}\times \Z$ subgroup of $H_1(\Sigma_\phi)$.

\begin{lem} \label{c1=0} 
Let $\xi:=\ker\alpha$. Then 
$c_1(\xi)=0$ on $H_2(M,\Z)$. Hence the $\mu$-index and 
$\bar{\mu}$-index for homologously trivial Reeb orbits 
are well defined.
\end{lem} 

\begin{proof} 
Apply Van Kampen theorem to $M=\Sigma_{\phi}\cup S^1\times D^2$ 
we have $H_2(M,\R)\cong H_2(\Sigma_\phi ,\R)$. This can be seen from 
the following long exact sequence: 
\begin{gather*} 
\cdots\to H_2(T^2)\overset{i_2}{\to}H_2(\Sigma_\phi)\oplus
H_2(S^1\times D^2)
\overset{j_2}\to H_2(M)\overset{\delta}{\to}  \\ 
H_1(T^2)\overset{i_1}\to H_1(\Sigma_\phi)\oplus H_1(S^1\times D^2)
\overset{j_1}{\to}H_1(M)\to \cdots
\end{gather*} 
Since $\text{im}(i_2)=0$ and $j_2=0$ on 
$H_2(S^1\times D^2)$, the inclusion 
$H_2(\Sigma_\phi)\hookrightarrow H_2(M)$ is injective. Now, 
$H_1(T^2)$ is generated by $[t]$ and $[B]$, and 
$i_1[t]=0$ in $H_1(S^1\times D^2)$, $i_1[t]\in H_1(\Sigma_\phi)$ 
is an element of infinite order; while $i_1[B]=0$ in 
$H_1(\Sigma_\phi)$, $i_1[B]\in H_1(S^1\times D^2)$ is a generator. 
So $i_1$ is injective, then $\text{im}(\delta)=0$, 
$\text{im}(j_2)=\ker(\delta)=H_2(M)$, and hence the inclusion 
$H_2(\Sigma_\phi)\hookrightarrow H_2(M)$ is surjective as well.
Hence  $H_2(M,\R )\cong H_2(\Sigma_\phi,\R)$ is 
generated by  
$\Gamma_\phi=\Gamma\times S^1_t$. Note that $\xi|_{\Gamma_\phi}$ 
is $S^1_t$-invariant, $\xi\cap T(\Gamma_\phi)$ generates a 
nonsingular 1-dimensional foliation, hence  
$c_1(\xi )[\Gamma_\phi]=0$ and we conclude that 
$c_1(\xi )=0$ on $H_2(M,\Z)$. 
\end{proof}

\begin{lem}
\begin{enumerate} 
\item Contractible Reeb orbits are 
\begin{itemize} 
\item $e_{k\sigma /m}$ and 
$h_{k\sigma/m}$ with $k,m\in\N$, $0<k<m$; 
\item $h^m$, $m\in \N$; and 
\item $B^m$, $m\in\N$, provided that $\sigma=1$ (note: $B^m$ is 
only homologously trivial but 
not contractible if $\sigma>1$). 
\end{itemize}   
\item If $n$ is not divisible by $\sigma$ then 
$[e_{n/m}]=[h_{n/m}]=(0,j)\in H_1(M,\Z)=\Z\oplus \Z_\sigma$ where 
$0\neq j\in\Z_{\sigma}$, $j\equiv n \mod \sigma$. 
\end{enumerate} 
\end{lem}

\vspace{.2in} 
\noindent
{\bf Moduli of planes and $\partial^2=0$.} \ 
Now the set of contractible Reeb orbits is identified. To 
determine whether or not the cylindrical contact homology 
is defined, i.e., whether or not $\partial^2=0$, we need 
to find all of the nonempty 
moduli $\cM(\gamma)$ with $\gamma$ contractible 
and $\dim\cM (\gamma)=1$.

Recall that if $\gamma$ is homotopically trivial and if 
$c_1(\xi)=0$ then 
\begin{equation} 
\dim \cM (\gamma )=\bar{\mu}(\gamma). 
\end{equation}  
So $\dim \cM (\gamma)=1$ implies that $\gamma$ is hyperbolic. 
Since $B^m$ is elliptic for all $m\in \N$ we have 
$\bar{\mu}(B^m)\equiv\bar{\mu}(B^m,\Z_2)=0$ $\mod 2$ and 
hence the following 

\begin{lem} \label{Bm} 
$\dim \cM(B^m)=\text{even}$ if $\cM(B^m)\neq\emptyset$. 
\end{lem} 

\begin{rem} 
In fact, $\bar{\mu}(B^m)$ can be 
made arbitrarily large by taking the constant $c$ in Lemma 
\ref{lambda'} to be large enough. 
\end{rem}

\begin{lem} 
$\cM(e_{k\sigma/m})=\emptyset =\cM(h_{k\sigma/m})$ for 
$k,m\in \N$, $0<k<m$. 
\end{lem} 

\begin{proof} 
Let $\gamma=e_{k\sigma/m}$ or $h_{k\sigma/m}$ with 
$k,m\in\N$, $0<k<m$. Assume that $\cM(\gamma)\neq \emptyset$ 
and let $C_M$ be the image in $M$ of an element of 
$\cM(\gamma)$. Since $m=\text{wind}(\gamma, B)\geq 1$ so 
$C_M\cap B\neq\emptyset$, $C_M$ intersects positively with $B$ 
at every point of intersection. We write $C\cap B=\{ z_i\}_{i=1}^s$ 
where $z_1,...,z_s$ are s distinct points.  
Let $m_i$ denote the intersection multiplicity at $z_i$, then 
$m_i>0$ and $\sum_{i=1}^s m_i=m$. Let $U_B$ be a thin tubular neighborhood 
of $B$. We may assume that $C\pitchfork \partial U_B$ 
and $C\cap U_B$ consists of $s$ disjoint discs $D_i\ni z_i$. 
Let $C':=C_M\setminus U_B$. We have 
$\partial C'=\gamma \cup(\cup_{i=1}^s(-\partial D_i))$. 
Note that each for reach $i$, $[\partial D_i] 
=[t]^{m_i}\in\pi_1(\Sigma_{\phi})$, while $[\gamma]=[p]^{k\sigma}[t]^m
\in \pi_1(\Sigma_{\phi})$. 
By Propositions \ref{-htpy+} and \ref{pi1} 
the existence of $C_M$ and hence $C'$ 
implies that $[p]^{k\sigma}=0\in\pi_1(\Sigma_{\phi})$ which 
is impossible unless $m\sigma$ divides $k\sigma$. 
By assumption $0<k<m$, $m\sigma$ does not divide $k\sigma$. 
So $C'$ does not exist. $\cM(\gamma)$ is empty. 
\end{proof}

Now consider $h$ the simple Type $\fC$ Reeb orbit. 
We can find a spanning disc $D_h$ for $h$ so that 
$D_{h}\cap B$ is a point, $\xi_{D_h}$ has only one 
singularity and it is elliptic. A simple calculation will 
yield the following 

\begin{lem}  \label{hm} 
For $m\in \N$, 
$\dim \cM(h^m)=\bar{\mu}(h^m)=2m-1$ provided that $\cM(h^m)\neq
\emptyset$. 
\end{lem} 

So $\dim \cM(h^m)=1\Rightarrow m=1$. One can show that, 
as the $S^1$-invariant case considered in \cite{Y2}, 
$\cM(h)/\R$ is in one-one correspondence with the gradient trajectories
from the corresponding 
hyperbolic critical point $x_h$ of $K$ to $\partial \Sigma$. 
There are two such trajectories and they are counted with opposite
signs. 

Hence we have 

\begin{lem}  \label{dh=0} 
$d_0 h=0$. (See Remark \ref{d} for the definition of $d_0$.) 
\end{lem} 

\begin{rem} 
{\rm 
Recall that to show $\partial^2=0$ it suffices to prove that 
$d_0=0$. Our computation above implies that indeed $d_0=0$ and 
hence $\partial^2=0$. The cylindrical contact homology is then defined. 
This will be confirmed again in Section \ref{orien} 
(see Lemma \ref{partial} and Corollary \ref{partial^2}).  
} 
\end{rem} 

\subsection{Holomorphic cylinders in $\R\times M$} \label{holocyl}

To compute the cylindrical contact homology of $(M,\xi)$ we 
need to find in $\R\times M$ 
all 1-dimensional moduli $\cM(\gamma_-,\gamma_+)$ 
of holomorphic cylinders converging to Reeb orbit $\gamma_-$ at 
$-\infty$ and to Reeb orbit $\gamma_+$ at $\infty$. 
First we find all ordered pairs $(\gamma_-,\gamma_+)$ such that 
$\dim \cM(\gamma_-,\gamma_+)=1$ if $\cM(\gamma_-,\gamma_+)
\neq\emptyset$, then we count elements of 
$\cM(\gamma_-,\gamma_+)/\R$.

\vspace{.2in} 
\noindent 
{\bf Type $\fB$ orbits.} \ 
The homotopy classes of Reeb orbits imply the 
following 

\begin{lem} 
$\cM(\gamma, B^m)=\emptyset =\cM(B^m,\gamma)$ for all $m\in\N$ if 
$\gamma$ is a Type $\fT$ or $\fC$ orbit. 
\end{lem} 

\begin{proof} 
Let $\cM=\cM(\gamma, B^m)$ or $\cM(B^m,\gamma)$. Assume that 
$\cM\neq\emptyset$ and let $C_M$ be the image in $M$ of an 
element of $\cM$. $C_M$ is a homotopy between $\gamma$ and 
$B^m$. Since $\gamma$ and $B^m$ are not homotopic in 
$\Sigma_\phi$, $C_M$ has to intersect $B$. 
We may assume that the intersection is transversal. By 
holomorphicity $C_M$ intersects positively with $B$ at every 
point of $C_M\cap B$. The total multiplicity 
of the intersection is $m=\text{wind}(\gamma, B)$. 
Thus we have 
$[\gamma][t]^{-m}=[B]$ in $\pi_1(\Sigma_\phi)$. 
But $[\gamma][t]^{-m}$ is in the abelian 
subgroup of $\pi_1(\Sigma_\phi)$ 
generated by $b=[\Gamma]$ while $[B]=aba^{-1}b^{-1}$ is not in 
this abelian subgroup, a contradiction to the existence of $C_M$. 
So $\cM=\emptyset$. 
\end{proof} 

\begin{lem} \label{BnBm} 
$\dim \cM(B^n,B^m)=\bar{\mu}(B^m)-\bar{\mu}(B^n)$ 
is an even number provided that $\cM(B^n,B^m)\neq\emptyset$. 
\end{lem} 

Recall that 
$\bar{\mu}(B^m)$ is always an even number and 
can be made arbitrarily large by modifying 
$\alpha$. This property together with Lemma \ref{Bm} and 
Lemma \ref{BnBm} shows that 
$B^m$ has eventually no contribution to the 
boundary operator $\partial$ of the cylindrical 
contact homology and hence 
can be neglected. 
We will then focus 
on Type $\fT$ and $\fC$ Reeb orbits 
and assume that $\gamma_{\pm}\neq B^m$, $\forall m\in \N$.

\vspace{.2in} 
\noindent 
{\bf Going upstairs.} \ 
Let 
\[ 
C\in\cM(\gamma_-,\gamma_+)
\] 
denote any one of such cylinders. 
After a 0-surgery along $B$, $C$ will become a punctured 
holomorphic sphere $\hat{C}\subset \R\times \hat{M}$ 
with one positive puncture and $1+s$ negative 
punctures, $s\geq 0$.

Recall that a 0-surgery along $B$ turns the open book 
$M=(\Sigma,\tau)$ into 
the mapping torus $\hat{M}$ with $B$ replaced by the orbit 
$\hat{e}$ corresponding to the extra elliptic fixed point of 
$\hat{\phi}'$ on $D^2=T^2\setminus \Sigma$;  
$e_{n/m}$, $h_{n,m}$ and 
$h^m$ ($0<n<\sigma m$, $n,m\in \N$) are lifted to $\hat{M}$ and are now
denoted by 
$\hat{e}_{n/m}$, $\hat{h}_{n/m}$ and $\hat{h}^m$ to  
distinguish them from their copies in $M$.

Surely the new holomorphic curve $\hat{C}$ 
converges to $\hat{\gamma}_+$ at $\infty$ and its 
original negative puncture converges to $\hat{\gamma}_-$ at 
$-\infty$. The number $s$ of extra negative punctures 
is equal to the the number of geometric intersection points 
(not counting multiplicity) of $C$ with $\R\times B$. 
Each of these punctures converges to $\hat{e}^{m_i}$ at 
$-\infty$ for some $m_i\in\N$. In other words, 
\[ 
\hat{C}\in \cM(\{ \hat{\gamma}_-,\hat{e}^{m_1},...,\hat{e}^{m_s}\},
\hat{\gamma}_+). 
\] 
$m_i$ is the intersection 
multiplicity of $C$ and $\R\times B$ at the corresponding point.

\begin{defn}[winding number around $B$]  
{\rm Suppose that $\gamma=h_{n/m}$ or $e_{n/m}$ or $h^m$. We call  
the number $m$ the {\em winding number} of $\gamma$ around $B$, 
and denote it by $\text{wind}(\gamma,B)$. 
}
\end{defn} 

We have 
\begin{equation}  
\sum_{i=1}^sm_i=\text{wind}(\gamma_+,B)-\text{wind}(\gamma_-,B).  
\end{equation}  

Note that $\gamma_+$ and $\gamma_-$ are free homotopic in $M$.

\vspace{.2in}
\noindent
{\bf One dimensional moduli of cylinders.} \ 
By using the homotopic property of Reeb orbits we can obtain 
the following

\begin{prop} \label{moduli} 
Let $\gamma_{\pm}$ be distinct Reeb 
orbits of Type $\fT$ or Type $\fC$. 
Suppose that $\cM(\gamma_-,\gamma_+)\neq \emptyset$ and 
$\dim \cM(\gamma_-,\gamma_+)=1$ then $(\gamma_-,\gamma_+)$ 
must be one of the following ordered pairs: 
\begin{enumerate}  
\item $(e_{n/m},h_{n/m})$, $0<n<\sigma m$; 
\item $(h_{n/(m-1)},e_{n/m})$, $0<n<\sigma (m-1)$; 
\item $(h^{m-1},e_{(\sigma m-1)/m})$; 
\item $(h_{(n-\sigma)/(m-1)},e_{n,m})$, $\sigma<n<\sigma m$; 
\item $(h^{m-1},e_{\sigma/m})$. 
\end{enumerate}  
\end{prop}

\begin{cor} \label{0-1} 
$\cC^o_*(\alpha )=0$ for $*\leq 0$. 
\end{cor} 

\begin{proof}[Proof of Corollary \ref{0-1}] 

Recall that $\cC^o(\alpha)$ is generated by 
(1) $B^m$ with $m\in\N$ if $\sigma=1$ and (2) 
$h^m$ with $m\in\N$, (3) $e_{k\sigma/m}$ and $h_{n/m}$ with 
$k,m\in\N$ and $0<k<m$. We have  known that 
$\bar{\mu}(B^m)\gg 1$ and $\bar{\mu}(h^m)=2m-1\geq 1$ for $m\in\N$. 
Now, Part 1 of Proposition \ref{moduli} implies that 
\[ 
\bar{\mu}(h_{k\sigma/m})=1+\bar{\mu}(e_{k\sigma/m}). 
\] 
By applying 
Part 2-5 of Proposition \ref{moduli} several times we get 

\begin{equation} \label{esigmam} 
\bar{\mu}(e_{k\sigma/m})=1+\bar{\mu}(h^{m-1})=2m-2\geq 2 \quad 
\text { for } k,m\in \N,  0<k<m, 
\end{equation}  
and hence $\bar{\mu}(h_{k\sigma/m})=2m-1\geq 3$ 
for  $k,m\in\N$ and $0<k<m$.  
 
This proves that $\cC^o_*(\alpha )=0$ for $*\leq 0$. 
\end{proof}

The rest of this section up to Lemma \ref{THS} 
is devoted to the proof of Proposition \ref{moduli}.

\begin{lem}
Assume that $\gamma_{\pm}$ are distinct 
Reeb orbits of Type $\fT$ or Type 
$\fC$.  Suppose that $\gamma_-$ and $\gamma_+$ are free
homotopic in 
$\Sigma_\phi$ and the formal dimension of $\cM(\gamma_-,\gamma_+)$ 
equals 1, then $(\gamma_-,\gamma_+)=(e_{n/m},h_{n/m})$. 
\end{lem}

Suppose that $\gamma_-$ and $\gamma_-$ are not free homotopic in 
$\Sigma_\phi$ and $\cM(\gamma_-,\gamma_+)\neq \emptyset$. 
Let $C\in \cM(\gamma_-,\gamma_+)$ be a holomorphic cylinder. Then 
$C\cap(\R\times B)\neq\emptyset$. Since $C$ and $\R\times B$ are 
holomorphic of complement dimensions, they intersects positively 
at every point of intersection. Let $C_M$ denote the image of 
$C$ in $M$. \\ 

\noindent 
{\bf Claim:} $C_M\cap h=\emptyset$. 

\begin{proof}[Proof of Claim:]
Let $L:=\gamma_+\cup (-\gamma_-)$. $[L]=[h]=0\in H_1(M,\Z)$ hence 
the linking number $lk(L,h)=lk(h,L)\in\Z$ is defined. Let $S$ be 
any surface with boundary $\partial S=L$. Then $lk(L,h)=S\cdot h$  
is the algebraic number of intersection points of $S$ with $h$. 
$lk(L,h)$ is independent of the choice of the spanning surface 
$S$ of $L$ (sine $h$ is homologously trivial). In particular, 
$C_M\cdot h=lk(L,h)$. On the other hand, we can find an embedded 
spanning disc $D$ of $h$ such that $D\cap B$ is a single point and 
$D\cdot L=\emptyset$. So $C_M\cdot h=lk(h,L)=D\cdot L=0$. Hence 
$C\cdot (\R\times h)=C_M\cdot h=0$. 
$C$ and $\R\times h$ are holomorphic 
of complement dimensions hence intersects positively at each point 
of intersection, thus we must have $C\cap(\R\times h)=\emptyset
=C_M\cap h$. 
\end{proof}

\begin{lem} \label{dimM-+} 
Assume that $\gamma_{\pm}$ are 
Reeb orbits of Type $\fT$ or Type 
$\fC$. Then 
\[ 
\dim
\cM(\gamma_-,\gamma_+)=r+
2(\text{wind}(\gamma_+,B)-\text{wind}(\gamma_-,B)),  
\] 
where 
\begin{enumerate} 
\item $r=0$ if $\gamma_{\pm}$ are both elliptic or hyperbolic; 
\item $r=1$ if $\gamma_+$ is hyperbolic and $\gamma_-$ is elliptic; 
\item $r=-1$ if $\gamma_+$ is elliptic and $\gamma_-$ is hyperbolic. 
\end{enumerate} 
\end{lem} 

\begin{proof} 

Recall that 
$C_M\cap B$ is a finite set of points say $z_1,...,z_s$. 
Let $m_i$ denote the multiplicity of $z_i$, then 
$\sum_{i=1}^sm_i=\text{wind}(\gamma_+,B)-\text{wind}(\gamma_-,B)
\geq 0$. Let $U_B$ denote a tubular neighborhood of $B\subset M$ 
so that $\partial U_B\pitchfork C_M$ and $U_B\cap C_M$ is a 
union of $s$ disjoint discs. 

On $\hat{\Sigma}=T^2$ take a nonvanishing vector field $Z_o$ and 
perturb it slightly to get a new vector field $Z$ such that
$Z$ has an elliptic singularity  at 
the critical point $x_e\in Crit(\hat{K})$ 
corresponding to $\hat{e}$, a hyperbolic singularity at the saddle point
$x_h\in Crit(\hat{K})$ corresponding to $\hat{h}$, and $Z$ is 
nonvanishing elsewhere. Since the volume form of $\Sigma=\hat{\Sigma}
\setminus D^2$ is positive on $\xi$, $Z$ defines a symplectic 
trivialization $\Psi$ on $M\setminus (U_B\cup h)$. 
The $\bar{\mu}$-index of $\gamma=\gamma_{\pm}$  with respect 
to $\Psi$ is 
\begin{equation}  
\bar{\mu}_\Psi (\gamma)=\begin{cases} 
-1 \quad \text{ if $\gamma =h_{n/m}$ or $h^m$}, \\ 
-2 \quad \text{ if $\gamma=e_{n/m}$}.  \end{cases} 
\end{equation}  

Let $\Psi'$ be a symplectic trivialization of $U_B\cong 
B\times D^2$ which is invariant under rotations in $B$. 
Let $D_i\subset (C\cap U_B)$ denote the connected component  
containing $z_i$, $i=1,...,s$. Consider along $\partial D_i$ 
the 1-dimensional subbundle $\cL_i:=\xi\cap T_{D_i}C$. $\cL_i$ is 
homotopic to the 1-dimensional subbundle spanned by the 
radio vector field $\partial_r|_{\partial D_i}$. 
Under the trivialization $\Psi$ the bundle $\cL_i$ becomes a 
loop in the Lagrangian Grassmann $\Lambda (\R^2)$ with 
Maslov index $n_i=0$. While with respect to $\Psi'$ the 
corresponding loop has Maslov index $n'_i=2m_i$ which is 
twice of the multiplicity of $z_i$. Thus we have 
(provided that $\cM(\gamma_-,\gamma_+)\neq \emptyset$) 
\begin{equation} 
\begin{split} 
\dim \cM(\gamma_-,\gamma_+) &=\bar{\mu}_{\Psi}(\gamma_+) 
-\bar{\mu}_{\Psi}(\gamma_-)+\sum_{i=1}^s(n_i'-n_i)\\ 
 &=\bar{\mu}_{\Psi}(\gamma_+) 
-\bar{\mu}_{\Psi}(\gamma_-)+\sum_{i=1}^s2m_i \\ 
&=\bar{\mu}_{\Psi}(\gamma_+) 
-\bar{\mu}_{\Psi}(\gamma_-)+2(\text{wind}(\gamma_+,B)-\text{wind}
(\gamma_-,B)).  
\end{split} 
\end{equation}  

Note that 
$|\bar{\mu}_{\Psi}(\gamma_+) -\bar{\mu}_{\Psi}(\gamma_-)|\leq 1$ 
and $\text{wind}(\gamma_+)-\text{wind}(\gamma_-)\geq 0$.  
So if $\dim \cM(\gamma_-,\gamma_+)=1$ then one of the followings 
must be held true: 
\begin{enumerate} 
\item $\gamma_+$ is hyperbolic, $\gamma_-$ is elliptic 
and $\text{wind}(\gamma_+,B)=\text{wind}(\gamma_-,B)$. 
\item $\gamma_+$ is elliptic, $\gamma_-$ is hyperbolic 
and $\text{wind}(\gamma_+,B)=1+\text{wind}(\gamma_-,B)$. 
\end{enumerate} 
This settles the denominators (i.e, winding numbers 
around $B$) of $\gamma_{\pm}$.

It is easy to see that $\gamma_+\neq h^m$, $\forall m\in \N$. 
The difference of the 
numerators of $\gamma_{\pm}$ is divisible by $\sigma$ 
since $\gamma_+$ and $\gamma_-$ are free homotopic. 
Assume that $\gamma_+=e_{n/m}$ and $\gamma_-=h_{n'/(m-1)}$ 
with $n-n'=k\sigma$. Note that notationally $h_{(im\sigma)/m}=h^m$ 
for all $i$. Then 
$k$ is the algebraic number of the times that 
the image $C_M\subset M$ of $C\in \cM(\gamma_-, 
\gamma_+)$ crosses $\Gamma_\phi$. Since $C_M$ intersects $B$ only once, 
we have $k=0$ or $1$.  This completes the proof. 
\end{proof}

\vspace{.2in} 
\noindent
{\bf Counting cylinders (up to signs).} \ Let  
\[ 
\kappa_{n,m}:=\text{gcd}(n,m) 
\] 
denote the greatest common divisor of $n$ and $m$. Recall 
that $\kappa_\gamma$ denotes the multiplicity of the Reeb 
orbit $\gamma$, so 
\[ 
\kappa_\gamma=\kappa_{n,m} \quad \text{ if $\gamma=e_{n/m}$ or 
$h_{n/m}$.} 
\]

\begin{lem} \label{THS} 
Let $(\gamma_-,\gamma_+)=(h_{n'/(m-1)},e_{n/m})$, 
$0<n<\sigma m$, $0\leq n'$, $n'=n$ or $n-\sigma$, 
be equal to one of four types of the ordered pairs 
in Proposition \ref{moduli}. Here $h_{0/(m-1)}$ represents the orbit
$h^{m-1}$.  Then the moduli $\cM(\gamma_-,\gamma_+)/\R$ consists of 
\[ 
\frac{\left| n'm-n(m-1)\right|}{\kappa_{n',m-1}
\kappa_{n,m}} 
\] 
points. All 
the corresponding 
pseudoholomorphic cylinders are {\em immersed} surfaces in 
$\R\times M$. 
\end{lem} 

\begin{rem} 
{\rm 
It will be proved later that a coherent orientation can be 
defined so that elements of the same moduli $\cM(
h_{n'/(m-1)},e_{n/m})$ have the same orientation. 
}
\end{rem} 

The rest of this section 
will be occupied by the proof of 
Lemma \ref{THS}. We will need Taubes's \cite{Ta} description on 
trice-punctured spheres, Bourgeois's \cite{B} Morse-Bott version 
of curve counting and some arguments brought from 
Hutchings-Sullivan's paper \cite{HS}.

\vspace{.2in} 
\noindent
{\bf More surgery models of $\hat{M}=\hat{\Sigma}_{\hat{\phi}}$.} \ 
Below we describe a more symmetric construction of $\hat{M}$ 
to compare with the contact structure on $S^2\times S^1$ considered 
in \cite{Ta}.  

Recall that $\hat{\Sigma}=T^2$ and $\hat{\phi}\in\text{Symp}(T^2,
\hat{\omega})$. We can choose $\hat{\omega}$ so that 
\[ 
\hat{\omega}=dq\wedge dp
\]  
with respect to some suitable coordinates 
$(q,p)\in \R/\Z\times \R/\Z$ of $T^2$. 
Then $\hat{\phi}$ is isotopic to $\psi\in\text{Symp}(T^2,
\hat{\omega})$, where 
\[ 
 \psi:T^2\to T^2, 
\quad \psi(q,p):=(q,p-\sigma q). 
\] 

The 2-form $\hat{\omega}\in\Omega^2(T^2)$ canonically extends 
to a 2-form on $T^2_{\psi}$ which is also denoted by 
$\hat{\omega}$. Note that $d\alpha\in\Omega^2(\Sigma_\phi)$ 
extends to $\hat{\omega}+\eta\wedge dt$ on $T^2_{\psi}$ 
for some closed 1-form $\eta\in \Omega^1(T^2)$.

Here is another way to think of $\hat{M}=T^2_\psi$. 
Think of $S^1_q$ as the interval $[-1,0]$ with 
two endpoints identified. Recall that the $q=0$ level is the 
mapping torus $(\{ 0\}\times S^1_p)_\psi$. 
Let $V=(-1,0)_q\times S_p^1$. 
Let $Z:=[-1,0]\times S^1\times S^1$. 
Following \cite{HS} we have the identification:  
\begin{align*} 
\Psi: V_{\psi} & 
\to (-1,0)\times S^1\times S^1 =\overset{\circ}{Z},\\ 
  [(q,p),t] & \to (q, p-\sigma qt, t). 
\end{align*} 
Then 
\begin{equation} \label{Zg} 
T^2_\psi\overset{\Psi}{\approx} 
\frac{[-1,0]\times T^2}{(-1,g(p,t):=(p+\sigma t,t))\sim 
(0,(p,t))},  
\end{equation} 
and the equality 
\[ 
\hat{\omega} =\Psi^*\omega _Z,\quad \omega_Z:=dq\wedge (dp+\sigma qdt),  
\] 
extends over $T^2_\psi$.

\vspace{.2in} 
\noindent
{\bf Orbits on $Z$ and in $T^2_\psi$.} \ The vector field 
$X_Z:=\partial_t-\sigma q\partial_p$ 
generates the line field $\ker(\omega_Z)$.  
Orbits of $X_Z$ comes in $S^1$-families 
parametrized by 
\begin{equation} \label{nm'} 
\sS_Z:=\{ (n,m)\in \Z_{\geq 0}\times \N \mid n\leq \sigma m \ 
 \}.  
\end{equation} 

\begin{notn} 
{\rm 
Let $\hU_{n/m}$ denote the $S^1$-family of orbits in $Z$ 
with index $(n,m)$. 
} 
\end{notn} 
The following facts are obvious: 
\begin{enumerate} 
\item The total space of $\hU_{n/m}$ is a 
$\kappa_{n,m}$-cover of the torus 
$\{ q=\frac{n}{\sigma m}\}$. 
\item For
$\gamma\in \hU_{n/m}$, 
$[\gamma]=(n,m)\in H_1(S^1_p\times S^1_t,\Z)$. 

\item 
With the gluing map $g$ in (\ref{Zg}) 
the two sets $\hU_{\sigma m/m}$ 
and $\hU_{0/m}$ are identified in $T^2_\psi$, and $S^1$-families 
of orbits on $T^2_\psi$ are indexed by the set 
\begin{equation} 
\hat{\sS}:=\{ (n,m)\in\Z_{\geq 0}\times\N \mid   0<n<\sigma m\}. 
\end{equation} 
\end{enumerate} 
Compare $\hat{\cS}$ with $\sS$ in (\ref{nm}) and one finds that 
$\hat{\sS}$ and $\sS$ are equal except the extra elements 
$(0,m)\in\hat{\sS}$ whose corresponding $S^1$-families of orbits 
will deform to the pair of orbits $\hat{h}^m, \hat{e}^m$ 
under a perturbation.

\vspace{.2in} 
\noindent 
{\bf Trice-punctured spheres in $\R\times Z$ and in $\R\times 
T^2_\psi$.} \ 
Recall that in \cite{Ta} (Thm. A.2) Taubes classified moduli of
trice-punctured  pseudoholomorphic spheres in 
the symplectization $\R\times (S^2\times S^1)$ 
equipped with  an almost complex structure associated to a certain 
contact structure on $S^2\times S^1$. 

Our manifold $Z=([-1,0]\times S^1)\times S^1$ can 
be suitably identified 
with a subset of $S^2\times S^1$. Note that 
with a constant 
$Q>0$ large enough, the 1-form 
$\alpha_Z:=qdp+(Q+\frac{1}{2}\sigma q^2)dt\in\Omega^1(Z)$ is contact 
and $d\alpha_Z=\omega_Z$, with Reeb vector field 
$R_Z=(Q-\frac{1}{2}\sigma q^2)^{-1}X_Z$. Then with a suitable 
embedding of $Z$ into $S^2\times S^1$, Taubes's contact 1-form 
on $S^2\times S^1$ restricts to $\alpha_T:=e^{\rho(q)}\alpha_Z$ on $Z$
for some $\rho\in C^{\infty}(Z)$ depending only on $q$. 
Also, up to a perturbation, $\alpha_T$ and $\alpha_Z$ have 
the same set of $S^1$-families of Reeb orbits indexed by 
(\ref{nm'}), the ordering of the orbits are preserved, and 
both orbit sets are invariant under rotations generated by 
the group $S^1_p\times S^1_t$. 

An energy 
inequality by Hutchings-Sullivan (see Sec. 3.3 of \cite{HS}) 
confirms that trice-punctured spheres converging to orbits in 
$Z$ stay in $\R\times Z$. Thus Taubes's result applies to 
$(Z,\alpha_Z)$.

Consider the triples $(n,m),(n',m'),(n'',m'')\in\sS_Z$ with 
$(n,m),(n',m')$ linearly independent and 
\begin{equation} \label{'''} 
n=n'+n'', \quad m=m'+m''.
\end{equation}    
Also use the notations 
\[ 
\hU_+:=\hU_{n/m},\quad \hU_-:=\hU_{n'/m'}, 
\quad \hU_*:=\hU_{n''/m''}. 
\] 
Recall that 
\[ 
\cM(Z;\hU_*,\hU_-;\hU_+) 
\] 
denote the moduli space of trice-punctured spheres 
in $\R\times Z$ with 
one positive end converging to an element of $\hU_+$, 
one negative end converging to an element of $\hU_-$, 
and one negative end converging to an element of $\hU_*$.

\begin{prop}[\cite{HS} Sec. 3.3] \label{confined} 
The total space of $\cM(Z;\hU_*,\hU_-;\hU_+)/\R$ 
is contained in $Z$ and is  
between $\{ q=\frac{n'}{\sigma m'}\}$ 
and $\{ q=\frac{n''}{\sigma m''}\}$. 
\end{prop} 

\begin{rem} 
{\rm 
Proposition \ref{confined} ensures that 
$\cM(Z;\hU_*,\hU_-;\hU_+)$ is preserved 
when boundaries of $Z$ are glued by $g$ to get $T^2_\psi$. 
}
\end{rem}

\begin{prop}[\cite{Ta} Thm. A.2] \label{TA2} 
$\cM(Z;\hU_*,\hU_-;\hU_+)/\R$ is connected and 
diffeomorphic to $S^1\times S^1$. 
The rotation group $G:=S^1_p\times S^1_t$ 
induces a free action on 
$\cM(Z;\hU_*,\hU_-;\hU_+)/\R$ 
with the orbit space a single point. Every element of 
$\cM(Z;\hU_*,\hU_-;\hU_+)$ is an immersed curve. 
\end{prop}

\vspace{.2in} 
\noindent
{\bf Counting in the Morse-Bott way.} \ 
Upon a perturbation using a suitable Morse function $f$
(see \cite{B}), 
each $S^1$-family of orbits deforms to a pair of 
orbits: one elliptic and one hyperbolic. 

Let 
\[ 
\he_*,\ \he_-,\ \he_+, \quad \hh_*,\ \hh_-,\ \hh_+ 
\] 
denote the elliptic and hyperbolic elements of $\hU_*, 
\hU_-,\hU_+$ respectively. 

We are interested in counting elements of 
0-dimensional moduli of the type 
$\cM(Z;\he_*,\hh_-;\he_+)/\R$. 

\begin{prop}[\cite{B} Sec. 3.2, 3.3] 
Assume that $\cM(Z;\he_*,\hh_-;\he_+)/\R\neq\emptyset$ and 
$\dim\cM(Z;\he_*,\hh_-;\he_+)/\R=0$. Then 
%
%
$\cM(Z;\he_*,\hh_-;\he_+)/\R$ is equal to the 
%
%
fibered product 
\[ 
\Big( W^s(\he_*)\times W^s(\hh_-)
\times_{(\hU_*
\times \hU_-)} \cM(Z;
\hU_*,\hU_-;\hU_+)\times 
_{\hU_+}W^u(\he_+)\Big)/\R,  
\] 
where $W^s(\cdot)$, $W^s(\cdot)$ are 
respectively the corresponding 
stable and unstable submanifolds of the Morse function;   
and the maps involved in the fibered product are evaluation 
maps. 
\end{prop} 
Note that $W^u(\he_+)=\{ \he_+\}$, $W^s(\hh_-)=\{\hh_-\}$ and 
$W^s(\he_*)=\hU_*\setminus \{\hh_*\}$.

\begin{rem} 
{\rm 
As the perturbation caused by the Morse function $f$ is small, 
elements of $\cM(Z;\he_*,\hh_-;\he_+)$ are immersed curves, following 
Proposition \ref{TA2}. 
} 
\end{rem} 

\begin{rem} 
{\rm 
it will be shown in Section \ref{orien} that elements of the 
same moduli $\cM(Z;\he_*,\hh_-;\he_+)$ have the same sign, hence 
up to sign the algebraic number of $\cM(Z;\he_*,\hh_-;\he_+)/\R$ 
equals the geometric number of $\cM(Z;\he_*,\hh_-;\he_+)/\R$. 
} 
\end{rem} 

To us the relevant cases are 
\begin{equation} \label{*} 
*=\begin{cases} (0,m-1) & \text{ if } 
0<n\leq \sigma (m-1), \\ 
(\sigma,m-1) & \text{ if } 
\sigma\leq n<\sigma m.\end{cases} 
\end{equation}  
Then the corresponding 1-dimensional spaces are 
\begin{enumerate} 
\item $\cM(Z;\he_{0/(m-1)},\hh_{n/(m-1)};\he_{n/m})$ for 
$0<n\leq \sigma(m-1)$; 
\item $\cM(Z;\he_{\sigma/(m-1)},\hh_{(n-\sigma)/(m-1)};\he_{n/m})$ for 
$\sigma\leq n<\sigma m$. 
\end{enumerate}

Recall the group $G$ from Proposition \ref{TA2}. 
Apply the group $G$ to the 
evaluation map 
\[ 
ev:\cM(Z;\hU_*,\hU_-;\hU_+)/\R 
\to \hU_{*-+}:=\hU_*\times \hU_-\times
\hU_+. 
\] 
We find that the image class $[ev(\cM/\R)]\in H_2(\hU_{*-+},\Z)$
is  equal to 
\[ 
\pm a[\hU_+\times \hU_-] 
\pm b[\hU_+\times \hU_*] 
\pm c[\hU_-\times \hU_*], 
\] 
for some $a,b,c\in\N$. 
\begin{fact} 
For a general triple $(n,m)=(n',m')+
(n'',m'')$ (with $(n,m),(n',m')$ linearly independent) 
we have $a=|nm'-mn'|(\kappa_{n,m}\kappa_{n',m'})^{-1}$.  

In particular 
\[ 
a=\begin{cases} 
n/(\kappa_{n,m}\kappa_{n,m-1}) \quad 
& \text{ if } \  *=(0,m-1), \\ 
(\sigma m-n)/(\kappa_{n,m}\kappa_{n-\sigma,m-1}) \quad 
& \text{ if } \ *=(\sigma,m-1). 
\end{cases} 
\] 
\end{fact}
Let 
\[ 
\hU_{*-+}\overset{\pi_*}{\to}\hU_*, \quad 
\hU_{*-+}\overset{\pi_{-+}}{\to}\hU_-\times 
\hU_+
\] 
denote the obvious projections. 
A direct calculation yields the following 
\begin{fact} 
Assume that (\ref{*}) holds. Then 
\[ \pi_*\Big( ev(\cM(Z;\he_*,\hh_-;\he_+)/\R)\cap 
\pi_{-+}^{-1}(\hh_-,\he_+)\Big)
\]  
consists of a single element in $W^s(\he_*)$. 
\end{fact}

By Lemma \ref{confined} holomorphic curve in $\R\times Z$ 
descent to holomorphic curves in $\R\times T^2_\psi$. 
Note that  
\[ 
\hh_{0/m}=\hh_{\sigma m/m}\quad \text{ and } \quad 
\he_{0/m}=\he_{\sigma m/m} \quad \text{ in } \hat{M}=T^2_\psi.
\]  
By Gromov compactness as in \cite{Hu} (Sec. 9.4),  
if a further deformation is applied then 
the numbers $\#\cM(\he_*,\hh_-,\he_+)/\R$ can change only when 
broken curves appear during the deformation. To us the 
only relevant cases are  
\begin{enumerate} 
\item broken curve from $\he_{n/m}$ to $\hh_{n/(m-1)}\hh$ 
to $\hh_{n/(m-1)}\he$, 
where the latter pieces consists of the trivial cylinder 
over $h_{n/(m-1)}$ and a cylinder between $h$ 
and $e$; 
\item broken curve from $\he_{n/m}$ to $\hh_{(n-\sigma)/(m-1)}\hh$ to 
$\hh_{(n-\sigma)/(m-1)}\he$,  
where the latter pieces consists of the trivial cylinder 
over $\hh_{(n-\sigma)/(m-1)}$ and a cylinder between $\hh$ 
and $\he$. 
\end{enumerate} 
In both cases the broken curves appear in cancelling pairs  
(due to the fact that $\cM(\he,\hh)/\R$ consists of two points 
with opposite signs, see Lemma \ref{dh=0}), 
so they do not affect the algebraic numbers of 
$\cM(\hat{M};\he_*,\hh_-;\he_+)/\R=\cM(M;\hh_-,\he_+)/\R$. 
This completes the 
proof of Lemma \ref{THS}.

\subsection{An energy estimate} \label{energy-est} 

Now we are back to the contact 3-manifold $(M,\xi)$. In the 
following 
we will derive an energy estimate which will help describe the 
signs of coefficients of the boundary operator $\partial$ 
(see Lemma \ref{c-c+}  in Section \ref{orien}).  

Let $\gamma=\gamma_{n/m}$be an unperturbed 
Type $\fT$ Reeb orbit with $0<n<\sigma m$. 
Recall that $\alpha =qdp +K(q)dt$. 
Then $\gamma\subset T^2_{q_o}=\{q_o\}\times S^1_p\times S^1_t$ with
$q_o$  satisfying $-K_q(q_o)=\frac{n}{m}$. Note that $\alpha$ 
pullbacks to a {\em closed} 1-form on $T^2_q\cong T^2$ that depends 
only on $q$. Thus for any closed curve $\gamma'\subset T^2_q$, 
the integral $\int _{\gamma'}\alpha$ depends only on the homology 
class of $\gamma'$. Consider the following homotopy of curves 
including $\gamma$: 
\begin{gather*} 
\Xi :[q_+,q'_+]\times S^1 \to [q_+,q'_+]\times S^1_p\times S^1_t, \\ 
\Xi (q,\theta ) = (q, p_o+ n\theta , m\theta). 
\end{gather*}   
We have 
\begin{equation} 
A_\gamma (q):=\int _{\Xi (q, \cdot )}\alpha =nq+mK(q).  
\end{equation} 
Differentiating $A_\gamma$ with respect to $q$ we have 
\[ 
\frac{dA_\gamma}{dq}=n+mK_q(q), \quad \frac{d^2A_\gamma}{dq^2}
=K_{qq}(q)>0.  
\] 
Note that $K_q(q_o)=\frac{-n}{m}$ and $K$ is a strictly decreasing 
function of $q$. $\frac{dA_\gamma}{dq}<0$ when $q<q_o$, 
$\frac{dA_\gamma}{dq}(q_o)=0$, and $\frac{dA_\gamma}{dq}>0$ when 
$q>q_o$. $q_o$ is the unique minimum point of $A_\gamma$. 
We have the following 

\begin{lem} \label{action-est} 
Let $\gamma=\gamma_{n/m}$be an unperturbed 
Type $\fT$ Reeb orbit with $0<n<\sigma m$. $\gamma\subset T^2_{q_o}$. 
Then the action function $A_\gamma(q)$ is a strictly increasing 
function of $|q-q_o|$. 
\end{lem}

\subsection{Orientation}  \label{orien} 

Let $J$ be an $\alpha$-admissible almost complex structure on 
$W:=\R\times M$ and let 
$\tilde{u}\in\cM (\gamma_-,\gamma_+)$ be a pseudoholomorphic map. 
$\tilde{u}:\R\times S^1\to W$. 
Let $(s,\tau)\in \R\times S^1$, $s+i\tau$ be the complex 
coordinate.  
%
%
$\tilde{u}$ satisfies the d-bar equation 
\[ 
\bar{\partial}\tilde{u}:=\tilde{u}_s+J\tilde{u}_\tau.  
\]

Assume from now on that $\tilde{u}$ 
is an {\em immersion}. Then there are two ways of splitting 
$\tilde{u}^*TW$ into a direct sum of two trivial 
$J$-complex line bundles: 
\begin{align*}  
\tilde{u}^*TW&=\tilde{u}^*\underline{\C}\oplus \tilde{u}^*\xi \\ 
   &=T\oplus N,  
\end{align*} 
where $\underline{\C}=\text{Span}(\partial_t,R_\alpha)$, 
$T=\text{Span}(\partial_s,\partial_\tau) =
\tilde{u}^*\text{Span}(\tilde{u}_s,\tilde{u}_\tau)$, 
and $N$ is the pullback of the normal bundle via $\tilde{u}$. 
We can compactify $\R\times S^1$ by adding to it two circles 
$S^1_{\pm\infty}$ at 
infinities to get $[-\infty, \infty]\times S^1$. 
Also compactify 
$\R\times M$ to get $[-\infty ,\infty]\times M$. 
$J$ can be extended over $\{\pm\infty\}\times M$ 
smoothly ,and $\tilde{u}$ can be 
extended over $S^1_{\pm\infty}$ pseudoholomorphically. 
The linearization of $\bar\partial$ at $\tilde{u}$ is 
\begin{gather*}  
D=D_{\tilde{u}}:W^{1,p}(\R\times S^1,\tilde{u}^*TW)\to 
L^p(\R\times S^1,\tilde{u}^*TW), \\ 
D\eta =\eta_s+J\eta_t+(\nabla_\eta J)\tilde{u}_\tau. 
\end{gather*}  

Trivialize $\tilde{u}^*TW=T\oplus N$ by the frame $\{ 
\partial_s,\partial_\tau, \nu, J\nu \}$ with 
$\nu, J\nu\in N$. Then 
\[ 
D\eta =\eta_s+J\eta_t+S(s,\tau)\eta,  
\] 
where $S(s,\tau)\in \R^{4\times 4}$ is of the block form 
\[ 
\begin{bmatrix} O& S_1\\O&S_2\end{bmatrix} 
\] 
with each block a $2\times 2$ matrix. Note that 
\[ 
S_1=O \text{ and } S_2=-J\nabla R_\alpha \quad \text{ at
}s=\pm\infty. 
\] 
Since $S_2$ is symmetric at $s=\pm\infty$, modulo some 
compact perturbation, we may as well assume that 
$S_1=O$ and $S_2$ is symmetric for all $s$ and $\tau$.

Write $D=\frac{d}{ds}-L$ where $L=-J\frac{d}{d\tau}-S$. 
Think of $L$ as an family of self-adjoint operators 
from $W^{1,p}(S^1,\R^4)$ to $L^p(S^1,\R^4)$. We can write 
$L=L_T\oplus L_N$ according to the decomposition $T\oplus N$. 
Then $L_T(s,\cdot )=-i\frac{d}{d\tau}$ for all $s$, and 
$L_N(s,\cdot )=-i\frac{d}{d\tau} -S_2(s,\cdot )$ with $S_2$ 
symmetric for all $s$. 

Since $D_{\tilde{u}}$ is surjective, $\text{Ind}(D_{\tilde{u}})
=\dim\ker(D_{\tilde{u}})=\dim \cM(\gamma_-,\gamma_+)$, which is 
the number of positive eigenvalues of $L(-\infty,\cdot )$, counted with 
multiplicity, flow to negative eigenvalues of $L(\infty,\cdot )$. 
Note that $L_T$ is static, i.e., independent of $s$, so the focus is 
on $L_N$. For each $s$, $L_N(s,\cdot )$ is a perturbation of 
the operator $-i\frac{d}{d\tau}$ by a $2\times 2$ symmetric matrix 
$S_2(s,\cdot )$. With a trivialization $\{ \nu,J\nu\}$ of $N$ fixed, 
following \cite{HWZ2} 
one can define the {\em winding number of an eigenvector field 
of $L_N(\pm\infty, \cdot )$ along $\gamma_{\pm}$} (relative to 
a trivialization of the normal bundle). We have the following 
properties from \cite{HWZ2}. 
\begin{itemize} 
\item For each $s$ the 
multiplicity of an eigenvalue of $L_N(s,\cdot )$ is at most 2. 
\item For each $k\in \Z$ 
and for each $s$ there are precisely two 
(not necessarily distinct) 
eigenvalues $\lambda_k^-\leq \lambda_k^+$ with 
winding number $k$, and the map from $s$ to the closed interval 
$[\lambda_k^-,\lambda_k^+]$ is continuous for every $k$, 
\end{itemize} 
The second property above means 
that the ordering of eigenvalues with distinct winding numbers is 
preserved under perturbation, and 
eigenvalues with distinct winding numbers never meet. This property 
greatly restricts the behavior of the 
spectral flow of $s\to L_N(s,\cdot )$.

A {{\em coherent orientation} of the moduli space is a choice of an 
orientation of the determinant bundle 
\[ 
\det(D_{\tilde{u}})=(\Lambda^{\text{max}}\ker(D_{\tilde{u}}))
\otimes (\Lambda^{\text{max}}
\text{coker}(D_{\tilde{u}}))^*
\] 
for all $\tilde{u}\in \cM(\gamma_-,\gamma_+)$ and all $\gamma_{\pm}$, 
such that these orientations match well under gluing operation 
(see \cite{FH}\cite{EGH}\cite{BM}). In the context of contact 
homology theory and symplectic field theory, coherent orientations 
always exist and may not be unique. 

Recall that 
0-dimensional moduli spaces must be of the form $\cM(\gamma,\gamma)$ 
and each of such moduli space consists of a single element whose 
image is $\R\times \gamma$. For $\tilde{u}\in\cM(\gamma,\gamma)$, 
$\ker (D_{\tilde{u}})=\{ 0\}=\text{coker}(D_{\tilde{u}})$, 
$\det(D_{\tilde{u}})=1\otimes 1^*$ has a natural orientation.

Note that for each Reeb orbit $\gamma$, the almost complex structure 
on $\xi$ induces a complex orientation on $\Gamma(\gamma,\xi_\gamma)
\cong C^{\infty}(S^1,\R^2)$. When $\alpha$ is regular, 
$L_{\gamma}=-J(\frac{d}{d\tau}-J\nabla R_\alpha |_\gamma )$ 
has no eigenvalue equal to 0, $\Gamma(\gamma,\xi_\gamma)$ 
can be written as a direct sum of two infinite dimensional vector 
spaces 
\[ 
\Gamma(\gamma,\xi_\gamma)=V_-(\gamma)\oplus V_+(\gamma),  
\] 
where $V_-(\gamma)$ (resp. $V_+(\gamma)$) is spanned by all 
negative (resp. positive) eigenvector fields of $L_{\gamma}$.

Given $\tilde{u}\in\cM(\gamma_-,\gamma_+)$, then the restriction 
of $\ker D_{\tilde{u}}$ to $\gamma$ is a subspace $E$ of 
$V_+(\gamma_-)$, such that $V_-(\gamma)\oplus E$ is identified with 
$V_-(\gamma_+)$ via the spectral flow. Certainly $E$ is sent by the 
flow to the restriction of $\ker D_{\tilde{u}}$ to $\gamma_+$, a 
subspace of $V_-(\gamma_+)$.

We apply the above discussion to our case here. Recall 
Proposition \ref{moduli}, Lemma \ref{THS} and Lemma \ref{partial}. 
We need to better understand the signs 
of terms in $\partial e_{n/m}$ to determine the boundary operator 
$\partial$. 

Let $\gamma$ be $h^m$ or 
a Type $\fT$ Reeb orbit
of either the unperturbed case ($\gamma=\gamma_{n/m}$) or 
the perturbed case ($\gamma=e_{n/m}$ or
$h_{n/m}$). $\gamma$ is 
the $s^{th}$ iterate of a simple Reeb orbit denoted by $\gamma'$. 
A tubular neighborhood of $\gamma'$ is 
diffeomorphic to $S^1\times D^2$, with $\gamma'$ identified with 
$S^1\times \{ 0\}$. Here we identify $S^1$ with $\R/a' \Z$ 
where 
$a'=\cA (\gamma')=\int_{\gamma'}\alpha $ is the action of 
$\gamma'$. 
Let $J$ denote an admissible almost complex structure. We 
may assume that $J|_{\gamma'}=\begin{bmatrix}0&-1\\1&0
\end{bmatrix}$ with respect to the ordered basis 
$\{ v_1=\partial_q,v_2=\partial_p-\frac{q}{K}\partial_t\}$. \\

\noindent
{\bf (a) A degenerated case.} 

Consider the differential operator $L_0:=-J\frac{d}{d\tau}$. 
The spectrum of $L_o$ is a discrete set $\text{spec}(L_0)
=\{ \frac{2k\pi a'}{s}\mid k\in \Z\}$. For each $\lambda\in 
\text{spec}(L_0)$ the eigenspace $V_{\lambda}(L_0)$ is a two 
dimensional vector space with an almost complex structure 
(hence is oriented) induced by $J$. A winding number 
(relative to the trivialization $v_1,v_2$) around 
$\gamma$ is defined for each eigenvector field of $L_o$. Two 
eigenvector field have the same winding number iff they belong 
to the same eigenspace $V_{\lambda}(L_0)$. In particular, the 
0-eigenspace $V_0(L_0)$ is spanned by two 
linearly independent eigenvector fields 
(i.e., $v_1$ and $v_2$) of winding number 0. \\

\noindent 
{\bf (b) Another degenerated case.} 

Recall (\ref{DR}) and consider the operator 
\[ 
L:=-J(\frac{d}{d\tau}-DR) \quad \text{ with }
DR=\begin{bmatrix}  0&0\\-\epsilon &0\end{bmatrix}, \ 0<\epsilon\ll 1 
\] 
$L$ corresponds to the case where Type $\fT$ Reeb orbits appear 
in $S^1$-families. $L$ is a small deformation of $L_0$ in {\bf 
(a)}  above. 
We may think that $L$ is 
included in a smooth homotopy $L_t$, $t\in [0,1]$ so that 
$L_1=L$. 
The winding numbers of eigenvector fields around $\gamma$ are 
unchanged under this homotopy in the sense that two 
eigenvector fields of $L_t$ belong to distinct eigenvalues of 
$L_t$ if their winding numbers around $\gamma$ are distinct. 
Thus the orientation of $V_{\lambda}(L_0)$ is transported along 
the homotopy to induce an orientation of the 2-dimensional 
vector space spanned by linearly independent eigenvector 
fields of $L$ with the same winding number around 
$\gamma$. 

We focus on the eigenvector fields of $L$ wit winding 
number 0. An easy computation shows that $v_1=\partial_q$ and 
$v_2=\partial_t-\frac{q}{K}\partial_t$ are two eigenvector 
fields with winding number 0, $\text{Span}(v_1)=V_\epsilon(L)$, 
$\text{Span}(v_2)=V_0(L)$. Note that 
\[ V_\epsilon(L)\oplus V_0(L)  
\text{ is homotopic to }V_0(L_0). 
\] 
So up to homotopy 
the perturbation term $DR$ has the 
effect of splitting $V_0(L_0)$ into a direct sum of a 1-dimensional 
vector space with positive eigenvalue and a 1-dimensional vector 
space with eigenvalue 0. 

Recall that a further deformation will turn the $S^1$-families of Reeb 
orbits into pairs $e_{n/m},h_{n/m}$. 
From the discussion following (\ref{DR}) we see that 
$V_\epsilon(L)$ will deform to $V_{\epsilon'}(L_\gamma)$ for some
$\epsilon'>0$ for $\gamma=e_{n/m},h_{n/m}$ and $h^m$. On the other 
hand, 
$V_0(L)$ will deform to $V_{\delta_+}(L_\gamma)$ for some $\delta_+>0$
if $\gamma=e_{n/m}$; to $V_{\delta_-}(L_\gamma)$ for some $\delta_-<0$
if $\gamma=h_{n/m}$ of $h^m$. This is summarized in {\bf (c),(d)} 
below. \\ 

\noindent 
{\bf (c) The elliptic case.} 

If $\gamma$ is elliptic then $V_-(L_{\gamma})$ is homotopic to 
$V_-(L_0)$ (and hence $V_-(L)$), 
$V_+(L_{\gamma})$ is homotopic to $V_0(L_0)\oplus 
V_+(L_0)$. So $V_{\pm}(L_\gamma)$ are endowed with 
induced orientations as described earlier. \\

\noindent 
{\bf (d) The hyperbolic case.} 

Now suppose that $\gamma$ is hyperbolic. 
Then $V_0(L_0)$ is 
deformed to a direct sum of two 1-dimensional eigenspaces 
\[ 
V_{\epsilon_+}(L_{\gamma})=\text{Span}(\tilde{v}_1)  \text{ 
and } V_{\epsilon_-}(L_{\gamma})=\text{Span}(\tilde{v}_2)
\]  
of 
$L_{\gamma}$, 
where $\epsilon_+$ (resp. $\epsilon_-$) is a the smallest positive 
eigenvalue (resp. the largest negative eigenvalue) of $L_{\gamma}$;  
$\tilde{v}_k$ is homotopic to $v_k$ for $k=1,2$. 
Moreover, 
\[ 
V_-(L_{\gamma})=V_{<\epsilon_-}(L_{\gamma})\oplus 
V_{\epsilon_-}(L_{\gamma}), \quad
V_+(L_{\gamma})=V_{\epsilon_+}(L_{\gamma})\oplus 
V_{>\epsilon_+}(L_{\gamma}),  
\]  
where $V_{<\epsilon_-}(L_{\gamma})$ and $V_{>\epsilon_+}(L_{\gamma})$ 
are with orientation induced from those of $V_-(L_0)$ and 
$V_+(L_0)$ respectively 
via deformation as before. 
Though $v_1\wedge v_2\sim \tilde{v}_1\wedge \tilde{v}_2$ 
induces a complex orientation for 
$V_{\epsilon_+}(L_{\gamma})\oplus V_{\epsilon_-}(L_{\gamma})
\cong \C$, in order to determine 
an orientation for each of $V_{\pm}(L_{\gamma})$ we still need to 
choose an orientation for each of $V_{\epsilon_{\pm}}(L_{\gamma})$ 
so that the combined orientation is consistent with the complex 
orientation. There are two choices: 
\begin{enumerate}  
\item either $o(V_{\epsilon_+}(L_{\gamma})) 
=o(\tilde{v}_1)$ and $o(V_{\epsilon_-}(L_{\gamma}))
=o(\tilde{v}_2)$, or 
\item $o(V_{\epsilon_+}(L_{\gamma}))
=o(-\tilde{v}_1)$ and $o(V_{\epsilon_-}(L_{\gamma}))
=o(-\tilde{v}_2)$. 
\end{enumerate}

It can be seen later that 
the vanishing of $\partial^2$ and, up to 
an isomorphism the resulting cylindrical contact homology 
are independent
of the choices of $o(V_{\epsilon_{\pm}}(L_{\gamma}))$. Here we just 
make a choice for all $\gamma=h_{n/m},h^m$ but we do not specify the 
choices as it 
will not hinder the computation of $HC(M,\xi)$. 

Let $\partial^*$ denote the dual operator of $\partial$ the 
boundary operator of cylindrical contact complex so that 
$\langle \partial\gamma_+,\gamma_-\rangle =\langle \gamma_+,
\partial^*\gamma_-\rangle$.

\begin{lem} \label{c-c+} 
We have for $m\geq 2$ 
\[ 
\partial^*h^{m-1}=\frac{1}{\sigma}
(c_-e_{\sigma(m-1)/m}+c_+e_{\sigma/m}). 
\] 
For $m\geq 2$ and $0<n<(m-1)\sigma$, 
\[ 
\partial^*h_{n/(m-1)}=
c_-\cdot\frac{\kappa_{n,m-1}}{n}e_{n/m}+c_+\cdot 
\frac{\kappa_{n,m-1}}{\sigma m-n-\sigma}e_{(n+\sigma)/m}
\] 
where $c_{\pm}\in \{ -1,1\}$ and $c_-c_+=-1$.  
\end{lem} 

\begin{proof} 
The only part that needs to be verified is the equality 
$c_-c_+=-1$. We will prove the second formula. The 
proof for $h^{m-1}$ is similar (using $h^{m-1}=h_{0/(m-1)}
=h_{\sigma(m-1)/(m-1)}$) and will be omitted. 

Apply 
the action estimate 
Lemma \ref{action-est} form Section \ref{energy-est} 
to the unperturbed case at first, i.e., the case of 
which Type $\fT$ Reeb orbits come in $S^1$-families indexed by 
$(n,m)$. Let $\Upsilon_{n/m}$ denote the $S^1$-family of Reeb 
orbits indexed by $(n,m)$ and let $\cM(\Upsilon_{n/(m-1)},
\Upsilon_{n'/m})$, $n'=n$ or $n+\sigma$, denote the 
moduli of pseudoholomorphic cylinders that converge to some 
element $\gamma_{n/(m-1)}\in\Upsilon_{n/(m-1)}$ at $-\infty$, 
and to some
element $\gamma'_{n'/m}\in\Upsilon_{n'/m}$ at $\infty$. 

Let $C_M\subset$ be the image in $M$ of some 
$\tilde{u}=(a,u)\in\cM(\Upsilon_{n/(m-1)},\Upsilon_{n'/m})$. 

Assume at first that $n'=n$. 
Then $\frac{n}{m-1}>\frac{n'}{m}$, $(n',m)=(0,1)+(n,m-1)$. 
$C_M\cap \Gamma=\emptyset$. 
$C_M$ is tangent to the vector field $v_1=\partial_q$ along 
the boundary $\gamma_{n/(m-1)}$. View $C_M$ as a homotopy of 
$\gamma_{n/(m-1)}$. Since $C_M$ is the image of a holomorphic 
curve, the integral of $\alpha$ along 
$C_M\cap T^2_q$ is increasing as $q$ increases from 
$q_0$ with $\gamma_{n/(m-1)}\in T^2_{q_0}$ to $q_1$ with 
$\gamma'_{n/m}\in\Upsilon_{n/m}$. 
Hence $u_s(-\infty, \cdot )$ is a {\em positive} multiple of 
$v_1=\partial_q$. The last property remains true after perturbation. 
I.e., if $\tilde{u}=(a,u)\in\cM(h_{n/(m-1)},e_{n/m})$ then 
$u_s(-\infty, \cdot )$ is a positive multiple of 
$\tilde{v}_1$. 
The flow $u_s$, $s\in \R$ induces an orientation on
$V_-(e_{n/m})$ that is compatible with the orientation 
\begin{equation}  
o_-:=o(\tilde{v}_1)\oplus o(V_-(h_{n/(m-1)}). 
\end{equation}   

On the other hand, if $n'=n+\sigma$ then 
$\frac{n'}{m}>\frac{n}{m-1}$ (recall that $n<(m-1)\sigma$), 
$(n',m)=(\sigma,1)+(n,m-1)$. $C_M$ has to cross $\Gamma_\phi$. 
By applying argument similar to the one in the case $n'=n$ 
we find that  if $\tilde{u}=(a,u)\in\cM(h_{n/(m-1)},
e_{(n+\sigma)/m})$ then 
$u_s(-\infty, \cdot )$ is a {\em negative} multiple of 
$\tilde{v}_1$. 
The flow $u_s$, $s\in \R$ induces an orientation on
$V_-(e_{(n+\sigma)/m})$ that is compatible with the orientation 
\begin{equation}  
o_+:=o(-\tilde{v}_1)\oplus o(V_-(h_{n/(m-1)}). 
\end{equation}

It is easy to see that {\em exactly one} of the 
following two equalities holds 
true: 
\[ 
o_-=o(V_-(e_{n/m})); \quad o_+=o(V_-(e_{(n+\sigma)/m})).  
\] 
So we have $c_-c_+=-1$. 
\end{proof}

Lemma \ref{c-c+}, 
Proposition \ref{moduli}, Lemma \ref{THS} and (\ref{coeff}) 
lead to the following

\begin{lem} \label{partial} 
The boundary operator $\partial$ of the contact complex of $(M,\xi)$ 
satisfies the following equations. 
\begin{enumerate} 
\item $\partial h^m=0=\partial B^m$ for $m\in \N$; 
\item $\partial h_{n/m}=0$ for $n,m\in\N$, $0<n<\sigma m$; 
\item $\partial e_{n/1}=0$  for $n\in N$, $0<n<\sigma$; 
\item $\partial e_{n/m}=c_-\cdot\frac{n}{\kappa_{n,m-1}}h_{n/(m-1)}$ 
for $0<n<\sigma$ and $m\geq 2$;  
\item $\partial e_{n/m}=c_+\cdot\frac{\sigma m-n}{\kappa_{n-\sigma,m-1}}
h_{(n-\sigma)/(m-1)}$ for 
$\sigma (m-1) <n<\sigma m$ and $m\geq 2$;  
\item $\partial e_{n/m}=c_-\cdot\frac{n}{\kappa_{n,m-1}} h_{n/(m-1)}
+c_+\cdot 
\frac{\sigma m-n}{\kappa_{n-\sigma,m-1}}h_{(n-\sigma)/(m-1)}$ 
for $\sigma\leq n\leq \sigma (m-1)$ and $m\geq 2$, 
\end{enumerate} 
where $c_-,c_+\in \{ -1,1\}$ and $c_-c_+=-1$. 
\end{lem} 

Note that the equation $\partial h_{n/m}=0$ can be obtained by 
a Mores-Bott argument. 

\begin{cor} \label{partial^2} 
$\partial ^2=0$. Hence the cylindrical contact homology 
$HC(M,\xi)$ is defined. Moreover, by Corollary \ref{0-1} 
and Theorem \ref{inv}, $HC(M,\xi)$ is independent of the 
choice of $(\alpha, J)$, and is an invariant of the 
isotopy  class of $\xi$. 
\end{cor}

\subsection{Computing $HC(M,\xi)$} \label{hc} 

Let $\cH$ denote the subcomplex generated by all hyperbolic 
Reeb orbits, and $\cE$ the subcomplex generated by all 
elliptic Reeb orbits (including no $B^m$). 

The  lemma below follows from Lemma \ref{partial}.  

\begin{lem} 
$\cH \subset\ker \partial$. In particular  
every hyperbolic Reeb orbit is a generator of $HC(M,\xi)$. 
Also, $\cE \cap \text{im}\partial =0$. 
\end{lem}

\begin{lem} 
$[h_{n/m}]=0\in HC(M,\xi)$ for all $0<n<\sigma m$ with 
$n$ not divisible by $\sigma$. 
\end{lem} 

\begin{proof} 
Let $0<i<\sigma$. 
From Lemma \ref{partial} we know that for $k=2,...,m-1$, 
$h_{(k\sigma-i)/m}$  is homologous 
(as an element of $HC(M,\xi)$) to a nonzero 
rational multiple of $h_{(\sigma-i)/m}$,  
and the latter is a nonzero rational multiple 
of $\partial e_{(\sigma-i)/m}$. This completes the proof. 
\end{proof}

\begin{lem} 
For $m\geq 1$ and $0<k<m$, 
$[h_{k\sigma/m}]$ is a nonzero rational multiple of $[h^m]$ 
in $HC(M,\xi)$. 
Also, $0\neq [h]\in HC(M,\xi)$. 
\end{lem} 

\begin{proof} 
The first statement follows from Lemma \ref{partial}. Use the 
equations  
\[ 
[h_{(k-1)\sigma/m}]=\frac{k\kappa_{(k-1)\sigma,m}}
{(m+1-k)\kappa_{k\sigma,m}}[h_{k\sigma/m}], \quad k=1,2,...,m, 
\] 
we get 
\[ 
[h^m] =[h_{0/m}] =\Big(\prod_{k=1}^{m}\frac{k}{m+1-k}\cdot 
\prod_{k=1}^m
\frac{\kappa_{(k-1)\sigma,m}}{\kappa_{k\sigma,m}}\Big)[h_{m\sigma/m}]
\] 
\[ 
=1\cdot 1 \cdot [h^m] 
\]  
There are no other boundary relations about $h^m$, so 
$0\neq [h]\in HC(M,\xi)$. 
\end{proof}

\begin{defn} \label{E_i} 
For $m\in\N_{\geq 2}$ define 
\[ 
E_{0,m}:=e_{\sigma/m}+\sum_{k=2}^{m-1}\Big( \prod_{j=1}^{k-1}
\frac{j}{(m-j-1)}\Big)e_{k\sigma/m}. 
\]    
Also for $i,m\in\N$ with $0<i<\sigma$ define 
\[   
E_{i,m}:=
\begin{cases} 
e_{i/1} & \text{ if } m=1, \\ 
e_{(\sigma-i)/m}+\sum_{k=2}^m\Big( \prod_{j=1}^{k-1}
\frac{j\sigma -i}{(m-j-1)\sigma+i}\Big)e_{(k\sigma-i)/m} 
& \text{ if } m\geq 2. 
\end{cases} 
\] 
\end{defn}

\begin{lem} 
$\partial E_{i,m}=0$ for all $(i,m)\in\Z_\sigma\times\N\setminus \{
(0,1)\}$. 
\end{lem} 

\begin{proof} 
Apply Lemma \ref{partial}. 
\end{proof}

\begin{lem} 
Let $E'$ be a finite 
linear combination of $e_{n/m}$'s with $\Q$-coefficients. 
Suppose that $\partial E'=0$. Then $E$ is a 
finite linear combination of 
$E_{i,m}$'s. 
\end{lem} 

\begin{proof} 
Write 
\[ \hspace{1in} 
E'=\sum_{n,m} c_{n,m}e_{n/m}=\sum_{i,m} E'_{i,m}, \qquad 
(i,m)\in\Z_\sigma\times\N\setminus \{
(0,1)\},  
\] 
where $E'_{i,m}$ consists of all $e_{n/m}$ terms of $E'$ 
with $n\equiv i \mod \sigma$. 
It is easy to see that 
\[ 
\partial E'=0 \text{ iff } \partial E'_{i,m}=0 \quad 
\forall i,m. 
\] 

Note that when $m=1$ and $0<i<\sigma$, each $E'_{i,1}$ is a rational 
multiple of $e_{i/1}=E_{i/1}$ 
(hence satisfies $\partial E'_{i,1}=0$), so we only need to consider 
the case $m\geq 2$. 
 
Fix $i,m$ with $m\geq 2$ and assume that $\partial E'_{i,m}=0$. 
Let $n_o$ 
be the smallest integer such that $n_o\equiv i\mod \sigma$
and $c_{n_o,m}\neq 0$. 

\noindent 
{\bf Claim:} $n_o\leq \sigma$. 

First assume instead that $n_o>\sigma (m-1)$, then 
$\partial E'_{i,m}=\partial (c_{n_o/m}e_{n_o/m})\neq 0$ 
by Lemma \ref{partial}, contradicting with $\partial E'_{i,m}=0$. 
So $n_o\leq \sigma(m-1)$.

Now suppose that $\sigma<n_o\leq \sigma(m-1)$, then 
\begin{align*} 
0 &=\partial E'_{i,m}  =\partial c_{n_o,m}e_{n_o/m}+
\partial \Big( \sum_{k=1}^{\lfloor m-\frac{n_o}{\sigma}\rfloor}
c_{(n_o+k\sigma),m}e_{(n_o+k\sigma)/m}\Big)  \\ 
& =\Big( c_0h_{(n_o-\sigma)/(m-1)}+c_1h_{n_o/(m-1)}\Big) + 
\Big( \text{no $h_{(n_o-\sigma)/(m-1)}$ term here }\Big) 
\\ 
& \neq 0 \quad \text{ because $c_0\neq 0$}. \quad \text{A 
contradiction again!} 
\end{align*}  
So $n_o\leq \sigma$.  

Suppose that $n_o<\sigma$, then $\partial e_{n_o/m}$ 
is a nonzero rational multiple of $h_{n_o/(m-1)}$ 
which can be cancelled only by adding a certain 
rational multiple of $\partial e_{(n_o+\sigma)/m}$, which 
again generates a nonzero rational multiple of $h_{(n_o+\sigma)/m}$ 
that can be cancelled only by adding a certain rational 
multiple of $\partial e_{(n_o+2\sigma)/m}$, and so on so forth. 
This process actually yields $E_{m,i}$ as defined 
in Definition \ref{E_i}. So $E'_{i,m}=c_{i,m}E_{i,m}$ for $i=1,...,
\sigma-1$. 
Similar arguments also apply to the case 
$n_o=\sigma$ and we have $E'_{0,m}=c_{\sigma,m}E_{0,m}$. 
This completes the proof. 
\end{proof} 

We have the following 

\begin{lem} \label{HC(M)} 
$HC(M,\xi)$ is freely generated by $[h^m]$, $m\geq 1$, and 
by $E_{i,m}$, $(i,m)\in\Z_\sigma\times\N\setminus \{ (0,1)\}$.  
The homologously trivial (as elements of $H_1(M,\Z)$) 
 generators are $[h^m]$, $m\geq 1$, and 
$E_{0,m}$ with $m\geq 2$. 
\end{lem} 

Note that $[h^m]=0=[E_{0,m}]\in H_1(M,\Z)$ 
(they are contractible actually), so their 
$\bar{\mu}$-indexes are defined. Recall from 
Lemma \ref{hm} that 
$\bar{\mu}(h^m)=2m-1$ for $m\in \N$. Also by (\ref{esigmam}) we have 
$\bar{\mu}(E_{0,m})=\bar{\mu}(e_{\sigma/m})=2m-2$ for 
$m\in\N_{\geq 2}$. 
This completes the proof of Theorem \ref{example}.

\end{document}